\documentstyle[12pt]{article}

\newtheorem{prop}{Proposition}[subsection]
\newtheorem{dfn}[prop]{Definition}
\newtheorem{theo}[prop]{Theorem}

\newtheorem{conj}[prop]{Conjecture}
\newtheorem{rem}[prop]{Remark}
\newtheorem{coro}[prop]{Corollary}
\newtheorem{exam}[prop]{Example}

\newtheorem{assume}[prop]{Assumption}

\def\C{{\bf C }}
\def\R{{\bf R }}

\def\H{{ \rm H }}
\def\Q{{ \bf Q }}
\def\P{{ \bf P }}
\def\N{{ \bf N }}
\def\NS{{ \rm NS}}

\def\a{ \alpha }
\def\b{ \beta }
\def\d{ \delta }
\def\e{ \varepsilon }
\def\g{ \gamma}

\def\k{ \kappa}

\def\L{ \Lambda }

\def\S{ \Sigma }

\def\ra{\rightarrow}

\title{\Huge Tamagawa numbers of polarized algebraic 
varieties }

\author{\small  
{\em Dedicated to Professor Yu.I. Manin on his 60th birthday} \\ 
$\; $\\
\sc Victor V. Batyrev and Yuri Tschinkel \\
$\; $ \\ 
$\; $\\
$\; $\\
}

\begin{document}

\date{}

\maketitle

\thispagestyle{empty}

\begin{abstract}
Let ${\cal L} = (L, \| \cdot \|_v)$ be an ample metrized invertible sheaf 
on a smooth quasi-projective algebraic variety $V$ over a number field $F$. 
Denote by $N(V,{\cal L},B)$ the number of rational 
points in $V$ having  ${\cal L}$-height $\leq B$.  
In this paper we consider  the problem of a  geometric 
and arithmetic interpretation of the asymptotic for  
$N(V,{\cal L},B)$   as $B \rightarrow \infty$ in connection 
with recent conjectures of Fujita  concerning the 
Minimal Model Program for polarized algebraic varieties.  

We introduce the notions  of ${\cal L}$-{\em primitive varieties} and
${\cal L}$-{\em primitive fibrations}.  For ${\cal L}$-primitive 
varieties $V$ over $F$ we propose a method to define an adelic 
Tamagawa number $\tau_{\cal L}(V)$ which is a generalization
of the Tamagawa number $\tau(V)$ introduced by Peyre for 
smooth Fano varieties. Our method allows us to  construct    
Tamagawa numbers for ${\bf Q}$-Fano varieties 
with at worst canonical singularities. 

In a series of examples of smooth polarized 
varieties and singular Fano varieties  we show that our 
Tamagawa numbers express the dependence  
of the asymptotic of $N(V,{\cal L},B)$  on the choice of 
$v$-adic metrics on ${\cal L}$.

\end{abstract}

\newpage

\tableofcontents

\newpage

\bigskip
\bigskip
\bigskip

\vskip 0,5cm

\section{Introduction} 
\bigskip 

Let $F$ be a number field (a finite extension 
of $\Q$),  ${\rm Val}(F)$  the set of all valuations of $F$,  
$F_v$ the $v$-adic completion of $F$ 
with respect to $v \in {\rm Val}(F)$, and  
$|\cdot |_v\, : \, F_v \ra \R$ the $v$-adic norm on $F_v$ normalized by 
the conditions $| x|_v = |N_{F_v/{\bf Q}_p}(x)|_p$ for $p$-adic 
valuations $v \in {\rm Val}(F)$.

Consider a projective space $\P^m$ 
with standard homogeneous 
coordinates $(z_0,...,z_m)$ and a locally closed 
quasi-projective 
subvariety $V \subset \P^m$ defined over $F$ (we want to stress that 
$V$ is not assumed to be projective).
Let $V(F)$ be the set of points in $V$ with 
coordinates in $F$. 
A  {\bf standard height function} $H\,:\, \P^m(F) \ra \R_{>0}$
is defined as follows
$$
H(x):=\prod_{v \in {\rm Val}(F)} \max_{j=0, \ldots, m} \{ |z_j(x)|_v   \}.
$$ 
A basic fact about the standard height function $H$ claims that 
the set 
\[ \{x\in \P^m(F)\,:\, H(x)\le B\} \] 
is finite for any real number $B$  \cite{lang}.  
We set   
$$
N(V,B)=\#\{x\in V(F)\; : \; H(x)\le B\}.
$$

It is an  experimental fact that whenever  one succeeds 
in proving an asymptotic formula for the 
function $N(V,B)$ as 
$B\ra \infty$,  one obtains the asymptotic 

\begin{equation}
N(V,B)= c(V) B^{a(V)}(\log B)^{b(V)-1}(1+o(1))
\label{formula}
\end{equation}
with some  constants $a(V) \in {\bf Q},$ $ b(V) 
\in \frac{1}{2}{\bf Z}$, and 
$c(V) \in {\bf R}_{>0}$. 
We want to use this observation as our 
starting point. It seems natural 
to ask the following: 
\medskip

{\bf Question A.} {\em For which 
quasi-projective subvarieties 
$V \subset {\bf P}^m$ defined over 
$F$ do there exist constants 
$a(V) \in {\bf Q},$ $ b(V) \in \frac{1}{2}{\bf Z}$  and 
$c(V) \in {\bf R}_{>0}$ such that the asymptotic 
formula $(1)$ holds? } 
\medskip

{\bf Question B.} {\em Does there exist a 
quasi-projective variety 
$V$ over $F$ with an asymptotic 
which is different from {\rm (\ref{formula})}?} 
\medskip 

In this paper we will be interested not in  
Questions $A$ and $B$ themselves but in a 
related to them  another  
natural question: 
\bigskip

\noindent
{\bf Question C.} {\em Assume that 
$V$ is an irreducible quasi-projective 
variety over a number field $F$ 
such that the   
asymptotic  formula {\rm (\ref{formula})} 
holds. How to compute the constants $a(V),b(V)$ and 
$c(V)$  in this formula via some  
arithmetical properties 
of $V$ over $F$ and geometrical properties 
of $V$ over ${\bf C}$? 
}
\bigskip

To simplify our terminology, it will be convenient for
us to postulate:
\bigskip

\noindent
{\bf Assumption.} For all quasi-projective $V', V$ with
$V'\subset V\subset \P^m$ and $|V(F)|=\infty $ 
there exists the limit
$$
\lim_{B\ra \infty} \frac{N(V',B)}{N(V,B)}.
$$
\bigskip

The following definitions have been useful to us:
\bigskip

\noindent 
{\bf Definition ${\bf S_1}$.} A smooth 
irreducible quasi-projective 
subvariety $V \subset {\bf P}^m$ over 
a number field $F$ 
is called {\bf weakly  saturated}, 
if $|V(F)| = \infty$ and if
for any locally closed 
subvariety $W \subset V$ with 
${\rm dim}\, W < {\rm dim}\, V$ one has 
\[ {\lim}_{B \rightarrow \infty} 
\frac{N(W,B)}{N(V,B)} < 1. \]
\bigskip

It is important to  
remark that  Question C really makes sense {\em only for 
weakly saturated} varieties. Indeed, if 
there were a locally closed 
subvariety $W \subset V$ with ${\rm dim}\, 
W < {\rm dim}\, V$ and  
\[ {\lim}_{B \rightarrow \infty} 
\frac{N(W,B)}{N(V,B)} = 1, \]
then it would be enough to answer Question 
C for each irreducible 
component of $W$ and for all possible  
intersections of these 
components (i.e., one could forget about 
the existence of $V$ and 
reduce the situation to a lower-dimensional case). 
In general, it is not 
easy to decide whether or not a  
given locally closed subvariety  
$V \subset {\bf P}^m$ is weakly saturated.  
We expect (and our assumption implies this) that
the orbits of connected subgroups 
$G \subset PGL(m+1)$ are
examples of weakly saturated varieties 
$V \subset {\bf P}^m$ (see \ref{equiv-sat}). 
\bigskip

\noindent
{\bf Definition} ${\bf S_2.}$ 
A smooth  irreducible quasi-projective subvariety 
$V \subset {\bf P}^m$ with $|N(V,B)| = \infty$ is 
called {\bf strongly saturated}, if 
for all dense Zariski open subsets 
$U \subset V$, one has 
\[  {\lim}_{B \rightarrow \infty}
 \frac{N(U,B)}{N(V,B)} = 1. \]
\bigskip

First of all, if $V \subset {\bf P}^m$ 
is a strongly saturated subvariety, 
then for any locally closed 
subvariety $W \subset V$ with ${\rm dim}\, 
W < {\rm dim}\, V$, one has 
\[ {\lim}_{B \rightarrow \infty} 
\frac{N(W,B)}{N(V,B)} =0, \]
i.e., $V$ is weakly saturated. 

On the other hand, if  $V \subset {\bf P}^m$ 
is weakly saturated, but not strongly saturated, 
then there must be an 
infinite sequence  $W_1, W_2, \ldots $ of 
pairwise different 
locally closed irreducible 
subvarieties  $W_i \subset V$ with ${\rm dim}\, 
W_i < {\rm dim}\, V$ 
and $|W_i(F)| = \infty$ 
such that for an arbitrary positive integer $k$ one has 
\[ 0 < {\lim}_{B \rightarrow \infty} 
\frac{N(W_1 \cup \cdots 
\cup W_k,B)}{N(V,B)} < 1. \]
Moreover, in this situation one can always 
choose the varieties 
$W_i$ to be strongly saturated 
(otherwise one could find $W_i' \subset W_i$ 
with  ${\rm dim}\, W_i' < {\rm dim}\, 
W_i$ with the same properties as 
 $W_i$ etc.). 
The strong saturatedness of each $W_i$ implies that 
\[ {\lim}_{B \rightarrow \infty} 
\frac{N(W_{i_1} \cap \cdots 
\cap W_{i_l},B)}{N(V,B)} = 0 \]
for all pairwise different $i_1, \ldots, i_l$ 
and $l \geq 2$. In particular, 
one has 
\[ \sum_{i=1}^k  {\lim}_{B \rightarrow \infty} 
\frac{N(W_i,B)}{N(V,B)} = 
 {\lim}_{B \rightarrow \infty} 
\frac{N(W_1 \cup \cdots \cup W_k,B)}{N(V,B)} 
< 1\;\; \forall k >0. \]
\medskip

\noindent
{\bf Definition} ${\bf F.}$ Let $V$ be a  
weakly saturated  quasi-projective 
variety in  ${\bf P}^m$ and $W_1, W_2, \ldots$ 
an infinite sequence 
of strongly saturated irreducible subvarieties 
$W_i$ having the property 
\[ 0 < \theta_i :=
 {\lim}_{B \rightarrow \infty} 
\frac{N(W_i,B)}{N(V,B)} < 1\;\; \forall i > 0. \]
We say that the set $\{W_1, W_2, \ldots  \}$
 forms an {\bf   
asymptotic arithmetic fibration} on $V$, 
 if the following equality holds
\[ \sum_{i=1}^{\infty} \theta_i = 1. \] 
\medskip

The main purpose of this paper is to explain some geometric and arithmetic 
ideas concerning weakly saturated varieties and their asymptotic 
arithmetic fibrations  by  strongly saturated subvarieties. 
It seems that the cubic 
bundles considered in \cite{BaTschi4} are 
examples of such a fibration. 
We want to remark that most of the   
above terminology grew out of our 
attempts to restore a conjectural picture of 
the interplay between the 
geometry of algebraic varieties and the 
arithmetic of the distribution 
of rational points on them 
after we have found in \cite{BaTschi4} 
an example which  contradicted  
general expectations formulated in \cite{BaMa}.  
\medskip

In section 2 we consider smooth quasi-projective 
varieties  $V$ over 
${\bf C}$ together with a polarization  
${\cal L} = (L, \| \cdot \|_h)$ consisting of 
an ample line bundle 
$L$ on $V$ equipped with a  positive 
hermitian metric $\| \cdot \|_h$.  
Our main interest in this section 
is a discussion of  geometric properties of 
$V$ in  connection  
with the Minimal Model Program \cite{KMM} and 
its version 
for polarized algebraic varieties suggested by Fujita 
\cite{fujita0,fujita01,fujita1}. We introduce 
our main geometric invariants $\alpha_{\cal L}(V)$, 
${\beta}_{\cal L}(V)$, and ${\delta}_{\cal L}(V)$ 
for an  arbitrary ${\cal L}$-polarized variety $V$. 
It is important to remark 
that  we will be  only interested in the case 
$\alpha_{\cal L}(V) > 0$. 
The number  $\alpha_{\cal L}(V)$ was first introduced in 
\cite{BaMa,Ba},  it equals to the opposite 
of the so called {\em Kodaira energy} 
(investigated by Fujita in
\cite{fujita0,fujita01,fujita1}). 
Our  basic geometric notion in the study 
of ${\cal L}$-polarized varieties $V$ with 
$\alpha_{\cal L}(V) >0$ is the notion 
of an ${\cal L}$-{\em primitive} variety. 
In  Fujita's program 
for polarized varieties with negative 
Kodaira energy ${\cal L}$-{primitive} 
varieties play the same role as  
${\bf Q}$-Fano varieties  in Mori's program  
for algebraic varieties with negative 
Kodaira dimension. In particular, 
one  expects the existence of so called 
${\cal L}$-{\em primitive fibrations}, 
which are analogous  
to  ${\bf Q}$-Fano fibrations in Mori's program.   
We show that 
on ${\cal L}$-primitive varieties there exists a 
canonical volume measure. Moreover, this 
measure allows us to construct 
a descent of hermitian metrics to the base of ${\cal L}$-primitive 
fibrations. Many  geometric ideas of  this section 
are inspired by \cite{BaMa,Ba}.   

In section 3 we introduce our main arithmetic notions of 
{\em weakly} and {\em strongly 
${\cal L}$-saturated} varieties. 
Our first main diophantine conjecture claims that 
if an adelic ${\cal L}$-polarized  
quasi-projective  algebraic variety $V$ 
over a number field $F$ is strongly 
${\cal L}$-saturated, then 
the corresponding ${\cal L}$-polarized 
complex algebraic variety  
$V({\bf C})$  is ${\cal L}$-primitive. 
Moreover, we conjecture  that  
if an adelic ${\cal L}$-polarized  
quasi-projective  algebraic variety $V$ 
over a number field $F$ is weakly  ${\cal L}$-saturated, then 
the corresponding ${\cal L}$-polarized complex algebraic variety  
$V({\bf C})$  admits an  ${\cal L}$-primitive fibration having infinitely 
many  fibers $W$ defined over $F$ which form an asymptotic arithmetic 
fibration.  
These conjectures allow us to establish 
a connection between the geometry of $V({\bf C})$ and the arithmetic 
of $V$.  Following this idea, we explain 
a construction of  an adelic measure on an arbitrary ${\cal L}$-primitive 
variety $V$ with ${\alpha}_{\cal L}(V) > 0$  
and of the  corresponding  Tamagawa number $\tau_{\cal L}(V)$  
as a regularized adelic integral of this measure.   
 Our construction generalizes  the definition 
of Tamagawa measures associated with a metrization of 
the canonical line bundle due to Peyre \cite{peyre}.   
We expect that for strongly 
${\cal L}$-saturated varieties $V$ 
the number $\tau_{\cal L}(V)$ reflects the dependence of 
the constant  $c(V)$ in the asymptotic formula (\ref{formula}) 
on the adelic metrization of the ample line bundle ${L}$. We discuss the
natural question about the behavior of the adelic constant 
$\tau_{\cal L}(W)$ for  fibers $W$ in   
${\cal L}$-primitive fibrations on weakly 
${\cal L}$-saturated varieties. 

In section 4 we show that our diophantine 
conjectures agree 
with already  known examples of asymptotic 
formulas established for polarized algebraic varieties 
through the study of analytic properties of 
height zeta functions.

In Section 5 we illustrate   our expectations 
for the constants $a(V), b(V)$ and $c(V)$ 
in the asymptotics 
of $N(V,B)$ on  some  examples of smooth Zariski 
dense subsets $V$ 
in  Fano varieties with singularities.

We would like to thank J.-L. Colliot-Th\'el\`ene
for his patience and encouragement. 
We are very grateful to  B. Mazur,  
Yu. I. Manin, L. Ein  and A. Chambert-Loir 
for their comments and suggestions. 
We thank the referee for several useful remarks.


\section{Geometry of ${\cal L}$-polarized varieties}

\subsection{${\cal L}$-closure }

Let $V$ be a smooth irreducible 
quasi-projective algebraic variety over ${\bf C}$, $V({\bf C})$ the 
set of closed points of $V$, 
$L$ an ample invertible sheaf on $V$,  i.e., $L^{\otimes k} 
= i^* {\cal O}_{{\bf P}^m}(1)$ for some $k >0$ and some  embedding 
$i \,: \, V \hookrightarrow {\bf P}^m$. Since we don't assume 
$V$ to be compact, the invertible sheaf $L$ on 
$V$ itself contains 
too little information about the embedding $i \,: \,
V \hookrightarrow {\bf P}^m$. For instance, 
let $V$ be an affine variety 
of positive dimension.  Then the space of 
global sections of $L$ 
is infinite dimensional and we don't know 
anything about the projective 
closure of $V$ in ${\bf P}^m$ even though we know that 
the invertible sheaf $L^{\otimes k}$ is isomorphic to 
$i^*{\cal O}_{{\bf P}^m}(1)$. This  situation changes if 
one considers $L$ together with a 
positive hermitian metric, i.e., 
an ample metrized invertible 
sheaf ${\cal L}$  associated with $L$. 
Let us choose a positive hermitian metric 
$h$ on ${\cal O}_{{\bf P}^m}(1)$ 
(e.g. Fubini-Study metric) and denote by $\|\cdot \|_h$ the induced 
metric on $L^{\otimes k}$. Thus  we obtain a metric $\| \cdot \|$ 
on $L$ by putting $\|s(x)\|: = \| s^k(x) \|_h^{1/k}$ for any $x \in 
V({\bf C})$ and any section 
$s \in \H^0(U, L)$ over an open subset $U \subset V$. 

\begin{dfn}
{\rm We call a pair ${\cal L} = ( L, \|\cdot \|)$ {\bf an ample  
metrized invertible sheaf}  associated with $L$.  We denote by 
${\cal L}^{\otimes \nu}$ the pair $( L^{\otimes \nu}, \|\cdot \|^{\nu })$.} 
\end{dfn}

Our next goal is to show that an ample metrized invertible sheaf 
contains almost complete information about the projective closure 
of $V$ in ${\bf P}^m$. 

\begin{dfn}
{\rm Let  ${\cal L} =( L,  \|\cdot \|) $ 
be an ample metrized invertible sheaf on a complex 
irreducible quasi-projective variety $V$. We denote by 
$$
\H^0_{\rm bd}(V, {\cal L})
$$ 
the subspace of $\H^0(V, L)$ 
consisting of those global sections $s$ of $L$ over $V$ 
such that the corresponding continuous 
function $x \mapsto \|s(x) \|$ $(x \in V({\bf C}))$  
is globally bounded on $V({\bf C})$ from above  by a  positive constant 
$C(s)$ depending only on 
$s$. We call $\H^0_{\rm bd}(V, {\cal L})$  the {\bf space of globally 
bounded sections} of ${\cal L}$. }
\end{dfn}

\begin{prop}
Let  $\overline{V}$ be the  normalization of the  projective closure
of $V$ with respect to the embedding  $i\; : \; 
V \hookrightarrow {\bf P}^m$ with  
${L}^{\otimes k} = i^*{\cal O}_{{\bf P}^m}(1)$. Denote by 
${e}\, : \,\overline{V} \ra {\bf P}^m$ the corresponding 
finite projective morphism.  Then 
one has  a natural isomorphism 
$$
\H^0_{\rm bd}(V, {\cal L}^{\otimes k})  \cong \H^0(\overline{V}, 
{e}^*{\cal O}_{{\bf P}^m}(1)).  
$$
\label{l-bd}
\end{prop}

\noindent
{\em Proof.} 
Since $\overline{V}({\bf C})$ is compact, the  
continuous function $x \mapsto \|s(x)\|$ is globally bounded on  
$\overline{V}({\bf C})$ for any 
$s \in \H^0(\overline{V}, {e}^*{\cal O}_{{\bf P}^m}(1))$. 
Therefore, we obtain that 
$\H^0(\overline{V}, {e}^*{\cal O}_{{\bf P}^m}(1))$ is a subspace 
of $\H^0_{\rm bd}(V, {\cal L}^{\otimes k})$. 

Now let $f \in 
\H^0_{\rm bd}(V, {\cal L}^{\otimes k})$ be a globally bounded on $V({\bf C})$ 
section of  $i^*{\cal O}_{{\bf P}^m}(1)$. Since 
$\H^0_{\rm bd}(V, {\cal L}^{\otimes k})$ is a subspace of 
$\H^0(V, L^{\otimes k})$,  the section $f$  uniquely 
extends to a global meromorphic section $\overline{f} 
\in \H^0(\overline{V}, {e}^*{\cal O}_{{\bf P}^m}(1))$. 
Since a bounded meromorphic function is holomorphic, 
$\overline{f}$ is a global regular section of 
${e}^*{\cal O}_{{\bf P}^m}(1)$ (we apply the theorem of Riemann to  
some resolution of singularities $\rho\, : \, X \ra \overline{V}$ 
and use the fact that
$\rho_* {\cal O}_X = {\cal O}_{\overline{V}}$). Thus we have 
$$
\H^0_{\rm bd}(V, {\cal L}^{\otimes k}) \subset  
\H^0(\overline{V}, {e}^*{\cal O}(1)). 
$$

\hfill $\Box$

\begin{dfn}
{\rm  We define a the graded ${\bf C}$-algebra 
$$ {\rm A}(V, {\cal L}) =  \bigoplus_{ \nu  \geq 0} 
\H^0_{\rm bd}(V, {\cal L}^{\otimes \nu } ).  \]
}
\end{dfn} 

Using \ref{l-bd}, one immediately obtains: 

\begin{coro}
The graded algebra  $ {\rm A}(V, {\cal L})$ 
is finitely generated. 
\end{coro}

\begin{dfn}
{\rm We call the normal projective variety 
\[ \overline{V}^{\cal L} = {\rm Proj}\,  
{\rm A}(V, {\cal L}). \]
the ${\cal L}$-{\bf closure} of $V$ with respect to 
an ample metrized invertible sheaf ${\cal L}$.
}
\end{dfn}

\begin{rem}
{\rm By \ref{l-bd}, $\overline{V}^{\cal L}$ 
is isomorphic to 
$\overline{V}$. Therefore,   
we have  obtained a way to define 
the normalization of the projective closure 
of $V$ with respect 
to an $L^{\otimes k}$-embedding via a notion 
of an ample  metrized invertible sheaf 
${\cal L}$ on $V$.}
\end{rem}

\subsection{Kodaira energy and $\alpha_{\cal L}(V)$}

Let $X$ be a normal irreducible  algebraic 
variety of dimension 
$n$.  We denote by ${\rm Div}(X)$ (resp. by 
${\rm Z}_{n-1}(X)$) 
the group of Cartier divisors 
(resp. Weil divisors) on $X$.
An element of ${\rm Div}(X) \otimes {\bf Q}$ 
(resp.${\rm Z}_{n-1}
(X) \otimes {\bf Q}$)  is called a 
${\bf Q}$-Cartier divisor (resp. a ${\bf Q}$-divisor). 
By $K_X$ we denote a  divisor of a meromorphic 
differential $n$-form 
on $X$, where  $K_X$ is considered as an 
element of  ${\rm Z}_{n-1}(X)$. 

\begin{dfn}
{\rm Let $X$ be a projective variety and 
$L$ be an invertible sheaf on 
$X$. The {\bf Iitaka-dimension} 
$\kappa(L)$ is defined as 
\[ \kappa(L) = \left\{ \begin{array}{ll} 
- \infty & \mbox{\rm if $\H^0(X, L^{\otimes \nu}) = 0$ for all $\nu > 0$} \\
{\rm Max}\, \dim \phi_{L^{\otimes \nu}}(X) : & 
 \H^0(X, L^{\otimes \nu})  \neq 0
\end{array} \right. 
\]
where $\phi_{L^{\otimes \nu}}(X)$ 
is the closure of the image of 
$X$ under the rational map 
\[ \phi_{L^{\otimes \nu}} \, : \, 
X \rightarrow {\bf P}(\H^0(X, L^{\otimes \nu})). \]
A Cartier divisor $L$ is called 
{\bf semi-ample} (resp. {\bf effective}), 
if $L^{\otimes \nu}$ is 
generated by global sections for some 
$\nu > 0$ (resp. $\kappa(L) \geq 0)$. 
}
\end{dfn} 

\begin{rem}
{\rm The notions of Iitaka-dimension,  
ampleness and semi-ample\-ness 
obviously extend to 
${\bf Q}$-Cartier divisors. 
Let $L$ be a Cartier divisor. Then 
for all  $\k_1, k_2 \in {\bf N}$ we set  
\[ \kappa(L^{\otimes k_1/k_2} ): = \kappa(L),\] 
\[   
 \mbox{\rm $L^{\otimes k_1/k_2}$ is  ample} \, \Leftrightarrow 
 \, \mbox{\rm $L$ is  ample}, \]
and 
\[    \mbox{\rm $L^{\otimes k_1/k_2}$ is  semi-ample} \, \Leftrightarrow 
 \, \mbox{\rm $L$ is  semi-ample}. \]
}
\end{rem}

\begin{dfn}
{\rm Let $X$ be a smooth projective variety. We denote by 
$\NS(X)$ the group of divisors on $X$ modulo 
numerical equivalence and set $\NS(X)_{\bf R} = \NS(X) \otimes {\bf R}$. 
By $[L]$ we denote the class of a divisor $L$ in $\NS(X)$. 
The {\bf cone of effective divisors} $\Lambda_{\rm eff}(X) \subset
\NS(X)_{\bf R}$ is defined as the closure of the subset 
\[ \bigcup_{ \kappa(L) \geq 0} {\bf R}_{\geq 0} [L] \subset 
\NS(X)_{\bf R}. \]
 }
\end{dfn}

\begin{dfn}
{\rm Let $V$ be a smooth quasi-projective algebraic variety 
with an ample  metrized invertible sheaf  ${\cal L}$, 
$\overline{V}^{\cal L}$ the ${\cal L}$-closure of $V$   
and $\rho$ some resolution of  
singularities 
\[ \rho \; : \; X  \rightarrow   \overline{V}^{\cal L}. \]
We define the number 
\[  \alpha_{\cal L}(V) =  \inf \{ t \in {\bf Q}\; : \; 
t [ \rho^*L] + [ K_X] \in  \Lambda_{\rm eff}(X) \}.  \] 
and call it the ${\cal L}$-{\bf index} of $V$. }
\end{dfn}

\begin{rem}
{\rm It is easy to see that the ${\cal L}$-{index} does not depend on 
the choice of $\rho$.} 
\end{rem}

\begin{rem}
{\rm The ${\cal L}$-index 
$\alpha_{\cal L}(V)$ for smooth projective 
varieties $V$ was first introduced in \cite{BaMa}  
and \cite{Ba}. We remark that  
the opposite number $- \alpha_{\cal L}(V)$ 
coincides with the notion 
of {\bf Kodaira energy} introduced and investigated by Fujita in 
\cite{fujita0,fujita01,fujita1}: 
\[ \kappa\epsilon (V, L) = - \alpha_{\cal L}(V) = 
 -  \inf \{ t \in {\bf Q}\; : \; 
\kappa ( (L)^{\otimes t} \otimes K_V) \geq 0 \}. \]
From the viewpoint  of our diophantine 
applications it is much more natural to consider  $\alpha_{\cal L}(V)$ 
instead of its opposite  $-\alpha_{\cal L}(V)$. The only reason that 
we could see for introducing  
the number $-\alpha_{\cal L}(V)$ instead 
of $\alpha_{\cal L}(V)$  is some kind of compatibility 
between the notions of {\em Kodaira energy} and {\em Kodaira dimension}, 
e.g. Kodaira energy must be positive (resp. negative) iff the 
Kodaira dimension is  positive (resp. negative). }
\end{rem} 

\noindent
The following statement was conjectured  in \cite{BaMa} (see also 
\cite{Ba,fujita00,fujita000}):

\begin{conj}
{\sc (Rationality)} Assume that  $\alpha_{\cal L}(V)> 0$. Then 
$\alpha_{\cal L}(V)$ is rational.
\end{conj}

\begin{rem}
{\rm 
It was shown in \cite{Ba} that this
conjecture follows from the Minimal Model Program. 
In particular, it holds for ${\rm dim}\,V \leq 3$. 
If ${\rm dim}\, V
=1$, then the only possible values of 
$\alpha_{\cal L}(V)$ are 
numbers $2/k$ with $(k \in {\bf N})$. If
 ${\rm dim}\, V =2$, then 
$\alpha_{\cal L}(V) \in \{ 2/k, 3/l \}$ 
with $( k, l \in {\bf N})$. }
\label{a-values}
\end{rem}

\begin{dfn}
{\rm A normal irreducible algebraic variety $W$ is said to have 
at worst {\bf canonical} (resp. {\bf terminal})  singularities 
if $K_W$ is a ${\bf Q}$-Cartier divisor 
and if for some (or every) resolution of singularities
\[  \rho \; : \; X  \rightarrow   W \] 
one has 
\[ K_X = \rho^*(K_W)\otimes {\cal O}(D) \]
where $D$ is an effective ${\bf Q}$-Cartier divisor (resp. 
the support of the effective divisor 
$D$ coincides with the exceptional locus of 
$\rho$). Irreducible components of the exceptional 
locus of $\rho$ which are not contained in the support of 
$D$ are called {\bf crepant divisors} of the resolution $\rho$.}
\end{dfn}

\begin{dfn}
{\rm  A normal irreducible algebraic variety $W$ is called 
a {\bf ca\-no\-nical} ${\bf Q}$-{\bf Fano} variety, if $W$ has  
at worst { canonical}   singularities 
and $K_W^{-1}$ is 
an ample ${\bf Q}$-Cartier divisor. A maximal positive 
rational number $r(W)$ such that $K_W^{-1} = L^{\otimes r(W)}$ 
for some Cartier divisor $L$ is called the {\bf index} of 
a canonical ${\bf Q}$-Fano variety $W$ (obviously, one has  
$r(W) = \alpha_{\cal L}(W)$ for some positive 
metric on $L$). }
\label{fano-c}
\end{dfn}

The following conjecture is due to Fujita \cite{fujita0}:

\begin{conj}
{\sc (Spectrum Conjecture)} Let $S(n)$ be the set all possible 
values of  $\alpha_{\cal L}(V)$ for smooth quasi-projective algebraic 
varieties $V$ of dimension $\leq n$ 
with an ample metrized invertible sheaf ${\cal L}$. Then 
for any $\varepsilon > 0$ the set 
\[ \{ \alpha_{\cal L}(V)  \in S(n) \; : \;  \alpha_{\cal L}(V) > 
\varepsilon  \} \]
is finite.
\end{conj} 

This conjecture follows from the Minimal Model Program \cite{KMM} and 
from the following conjecture on the 
boundedness of index for Fano varieties with canonical singularities:
 
\begin{conj}
{\sc (Boundedness of Index)} The set of possible values of index
$r(W)$ for canonical ${\bf Q}$-Fano varieties $W$ of dimension $n$ 
is finite. 
\end{conj}

In particular, both conjectures are true for 
${\bf Q}$-Fano varieties of dimension 
$n \leq 3$ \cite{A,fujita0,fujita01,fujita1,Ka,shin}.

\subsection{${\cal L}$-primitive varieties} 

\begin{dfn}
{\rm Let $X$ be a  projective algebraic variety.   
We call  an effective ${\bf Q}$-divisor $D$ {\bf rigid}, if 
$\kappa(D) = 0$. 
}
\end{dfn}

\begin{prop}
{\rm An effective  ${\bf Q}$-divisor 
$D$ on $X$ is  rigid if and only if there exist finitely many 
irreducible subvarieties  $D_1, \ldots, D_l \subset X$ $(l \geq 0)$ 
of codimension $1$ such that 
$D =r_1 D_1 +  \ldots + r_l D_l$ with $r_1, \ldots, r_l \in {\bf Q}_{>0}$ 
and 
\[ \mbox {\rm dim}\, \H^0 (X, 
{\cal O}(n_1 D_1 +  \ldots + n_l D_l)) =1 \;\; \forall\; 
(n_1, \ldots, n_l ) \in {\bf Z}^l_{\geq 0}.  \]
}
\label{rigid2}
\end{prop}

\noindent
{\em Proof.} Let $D$ be rigid. Take a positive integer $m_0$  
such that $m_0D$ is a Cartier divisor and 
${\rm dim}\, \H^0(X, {\cal O}(m_0D) ) =1$. Denote 
by $D_1, \ldots, D_l$ the irreducible components of the divisor $(s)$ of 
a non-zero section $s \in  \H^0(X, {\cal O}(mD) )$.  
One has 
\[ (s) = m_1 D_1 + \cdots + m_l D_l\;\;\; m_1, \ldots, m_l \in {\bf N}. \]
Since ${\cal O}(D_i)$ admits at least one global non-zero section we obtain that 
\[ \mbox {\rm dim}\, \H^0 (X, 
{\cal O}(n_1' D_1 +  \ldots + n_l' D_l)) \geq 
\mbox {\rm dim}\, 
\H^0 (X, {\cal O}(n_1 D_1 +  \ldots + n_l D_l)), \]
whenever $n_1 \geq n_1', \ldots, 
n_l \geq n_l'$ for $(n_1', \ldots, n_l' ), \; 
(n_1, \ldots, n_l )  \in {\bf Z}^l_{\geq 0}$. 
This implies that    
\[ \mbox {\rm dim}\, \H^0 (X, 
{\cal O}(n_1 D_1 +  \ldots + n_l D_l)) \geq 
1\;\;  \forall\; 
(n_1, \ldots, n_l ) \in {\bf Z}^l_{\geq 0}.  \]
On the other hand, for any  $(n_1, \ldots, n_l )  
\in {\bf Z}^l_{\geq 0}$
there exists a positive integer $n_0$ such that 
$n_0m_1 \geq n_1, \ldots, n_0m_l \geq n_l$. Therefore, 
\[ \mbox {\rm dim}\, \H^0 (X, 
{\cal O}(n_1 D_1 +  \ldots + n_l D_l)) \leq 
\mbox {\rm dim}\, \H^0 (X, 
{\cal O}(n_0m_0D) )  =1, \]
since $\kappa(n_0m_0D) =\kappa(D) =0$. 

\hfill $\Box$ 

\begin{coro}
Let   $D_1, \ldots, D_l \subset X$ be 
all irreducible components of the 
support of a rigid ${\bf  Q}$-Cartier divisor $D$. Then 
a linear combination 
\[ n_1 D_1 +  \ldots + n_l D_l, \; \; n_1, \ldots, n_l \in {\bf Z} \]
is a principal divisor, iff $n_1 = \cdots = n_l = 0$. 
\label{rigid-l}
\end{coro}

\noindent
{\em Proof.}
 Assume that $n_1 D_1 +  \ldots + n_l D_l$ is linearly 
equivalent to $0$. Then the effective Cartier divisor  
$D_0 = \sum_{n_i \geq 0} n_i D_i$ is linearly equivalent 
to the effective Cartier 
divisor $D_0' = \sum_{n_j <0} (-n_j) D_j$. 
Since $D_0$ and $D_0'$ have different supports we have 
${\rm dim}\, \H^0(X, {\cal O}(D_0)) \geq 2$. 
Contradiction to \ref{rigid2}. 

\hfill $\Box$

\begin{dfn}  
{\rm Let $V$ be a smooth quasi-projective algebraic variety 
with an ample  metrized invertible sheaf  ${\cal L}$ and 
 $\overline{V}^{\cal L}$ the  projective ${\cal L}$-closure of $V$. 
The variety $V$ is called  ${\cal L}$-{\bf primitive}, if  
the number $\alpha_{\cal L}(V)$ is rational and if for some resolution of  
singularities 
\[ \rho \; : \; X  \rightarrow   \overline{V}^{\cal L} \]
one has $\rho^*(L)^{\otimes \alpha_{\cal L}(V)} \otimes K_X = 
{\cal O}(D)$, where $D$  is a rigid  
effective ${\bf Q}$-Cartier divisor  on $X$. 
} 
\end{dfn}

\begin{rem}
{\rm It is easy to see that the notion  of an ${\cal L}$-primitive 
variety  doesn't depend 
on the choice 
of a resolution of singularities $\rho$.   Since $V$ is smooth,  
we can  always assume that 
the natural mapping 
$$\rho\; : \; \rho^{-1}(V) \rightarrow V$$
is an isomorphism. 
}
\end{rem}

\begin{exam}
{\rm Let  $V_1$ and $V_2$ be two smooth quasi-projective 
varieties with ample metrized invertible sheaves 
${\cal L}_1$ and ${\cal L}_2$ (resp. on $V_1$ and $V_2$).  
Assume that $V_1$ (resp. $V_2$ ) is ${\cal L}_1$-primitive 
(resp.  ${\cal L}_1$-primitive).  Then the product 
$V = V_1 \times V_2$ is ${\cal L}$-primitive, where   
${\cal L} = \pi_1^*{\cal L}_1 \otimes \pi_2^*{\cal L}_2$. }
\end{exam}

\noindent
Our main list of examples of ${\cal L}$-primitive varieties is obtained 
from canonical ${\bf Q}$-Fano varieties:

\begin{exam}
{\rm Let $V$ be the set of nonsingular 
points of a canonical 
${\bf Q}$-Fano variety $W$ with an ample  metrized 
invertible  sheaf ${\cal L} = (L,\|\cdot\|)$ 
such that  $K_W^{-1} = L^{\otimes r(W)}$  ( $W = 
\overline{V}^{\cal L} $). Then $V$ is 
an ${\cal L}$-primitive variety with 
${\cal L}$-index $r(W)$. Indeed, let $\rho\, : \, X \ra W$ be a resolution 
of singularities. By \ref{fano-c}, we have 
$\rho^*(L)^{\otimes r(W)} \otimes 
K_X = {\cal O}(D)$, where 
$D$ is an effective 
${\bf Q}$-Cartier divisor. Since  the 
support of $D$ consists of exceptional 
divisors with respect to $\rho$, $D$ is rigid 
(see \ref{rigid2}). }
\end{exam}

We expect that the above examples cover  
all  ${\cal L}$-primitive varieties: 

\begin{conj}
{\sc (Canonical ${\bf Q}$-Fano contraction)}  
Let $V$  be an ${\cal L}$-primitive variety with 
$\alpha_{\cal L}(V) >0$. Then there exists  
a resolution of singularities $\rho\, : \, X \ra \overline{V}^{\cal L}$ 
and a birational projective morphism $\pi\, : \, X \ra W$ to a 
canonical ${\bf Q}$-Fano variety $W$ such that 
$\pi^*K_W^{-1} \cong \rho^*(L)^{\alpha_{\cal L}(V)}$ 
(i.e., $\alpha_{\cal L}(V) = r(W)$) 
and the support of $D$ 
($\rho^*(L)^{\otimes r(W)} \otimes 
K_X = {\cal O}(D)$) is contained 
in the exceptional locus of $\pi$.
\label{conj-cont}
\end{conj}

The above conjecture is expected to follow 
from the Minimal Model Program 
using the existence and  termination of 
flips (in particular, it holds for toric varieties).  

The following statement will be important in
our construction of Tamagawa numbers 
for ${\cal L}$-primitive varieties 
defined over a number field: 

\begin{conj} {\sc (Vanishing)}  
{\rm For  $V$ an ${\cal L}$-primitive variety such that 
$\alpha_{\cal L}(V) >0$ we have 
\[ {\rm h}^i(X, {\cal O}_X) = 0 \;\; \forall \; i > 0 \]
for any resolution of singularities 
$$
\rho \; : \; X  \rightarrow   \overline{V}^{\cal L}
$$
such that the support of the ${\bf Q}$-Cartier divisor 
$\rho^*(L)^{\otimes \alpha_{\cal L}(V) } 
\otimes K_X$ is a ${\bf Q}$-Cartier divisor with normal crossings. 
In particular, $ {\rm Pic}(X)$ is a finitely generated 
abelian group and  one has a  canonical isomorphism 
\[ {\rm Pic}(X) \otimes {\bf Q} \cong \NS(X) \otimes {\bf Q}. \]
\label{vanish}
} 
\end{conj}

\begin{rem}
{\rm Theorem 1-2-5 in  \cite{KMM} implies the vanishing for 
the structure sheaf for ${\bf Q}$-Fano varieties with canonical 
singularities (even with log-terminal singularities). 
All canonical (and log-terminal singularities) are rational and it
follows that the higher cohomology of the structure sheaf 
on any desingularization of a canonical 
or a log-terminal Fano variety must also vanish 
(by Leray spectral sequence).
Therefore, we would obtain the vanishing \ref{vanish} for all 
${\cal L}$-primitive varieties
which are birationally equivalent to a 
Fano variety with at worst log-terminal 
singularities. The existence of a canonical ${\bf Q}$-contraction  
\ref{conj-cont} would insure this.
\label{vanish-rem}
}
\end{rem}

\begin{dfn}
{\rm Let $V$  be an ${\cal L}$-primitive variety with 
$\alpha_{\cal L}(V) >0$, $\rho \, : \, X  \rightarrow   
\overline{V}^{\cal L}$ any resolution of singularities, 
$D_1, \ldots, D_l$ irreducible components of the support of the 
rigid effective
${\bf Q}$-Cartier divisor $D$ with ${\cal O}(D) =
\rho^*(L)^{\otimes \alpha_L(V)} \otimes K_X$. 
We shall call 
$$
{\rm Pic}(V, {\cal L}) : = {\rm Pic}(X \setminus \bigcup_{i =1}^l D_i)
$$
the ${\cal L}$-{\bf Picard group} of 
$V$. The number  
\[ \beta_{\cal L}(V): = {\rm rk}\, {\rm Pic}(V, {\cal L})  \]
will be called the ${\cal L}$-{\bf rank} of $V$. 
We define the ${\cal L}$-{\bf cone of effective divisors} 
$\Lambda_{\rm eff}(V, {\cal L}) \subset
{\rm Pic}(V,{\cal L})\otimes {\bf R}$ as the image  
of $\Lambda_{\rm eff}(X) \subset \NS(X)_{\bf R} = 
{\rm Pic}(X) \otimes {\bf R}$ under the natural surjective 
${\bf R}$-linear mapping 
\[ \tilde{\rho}\; : \;{\rm Pic}(X) \otimes {\bf R} \ra 
{\rm Pic}(V,{\cal L})\otimes {\bf R}. \]
}
\end{dfn}

\begin{rem}
{\rm 
By \ref{rigid-l}, one obtains the  
exact sequence 
\begin{equation}
 0 \ra {\bf Z}[D_1] \oplus \cdots \oplus  
{\bf Z}[D_l] \ra {\rm Pic}(X) \stackrel{\tilde{\rho}}{\ra}
  {\rm Pic}(V, {\cal L}) {\ra} 0 
\end{equation}
and therefore 
\[  \beta_{\cal L}(V)= 
{\rm rk}\, {\rm Pic} (X) - l.  \]
Using these facts, it is easy to show that  the group 
$ {\rm Pic}(V, {\cal L})$ and the cone 
$\Lambda_{\rm eff}(V, {\cal L})$ do 
not depend on the choice of 
a  resolution of singularities 
$\rho \, : \,X  \rightarrow   \overline{V}^{\cal L}$. }
\end{rem}

\noindent
The above conjecture holds in 
dimension $n \leq 3$ as a consequence 
of the Minimal Model Program. More precisely, 
it is a consequence 
of Conjecture \ref{conj-cont} and the 
following weaker statement:

\begin{conj} 
{\sc (Polyhedrality)} 
Let $V$  be an ${\cal L}$-primitive variety with 
$\alpha_{\cal L}(V) >0$. Then 
$\Lambda_{\rm eff}(V, {\cal L})$
is a rational finitely generated polyhedral cone. 
\label{polyhed} 
\end{conj} 
\medskip

\begin{dfn}{\rm 
Let $ (A,A_{\R},\L)$ be a triple consisting
 of a finitely generated 
abelian group
$A$ of rank $k$, a $k$-dimensional real vector space $A_{\R}=A\otimes \R$
and a convex $k$-dimensional finitely generated  polyhedral cone
$\L\in A_{\R}$ such that $\L\cap -\L=0\in A_{\R}$. 
For ${\rm Re}({\bf s})$ contained in the interior of the cone $\L$
we define the ${\cal X}$-function of $\L$ by 
the integral
$$
{\cal X}_{\L}({\bf s}):= \int_{\L^*}e^{-<{\bf s},{\bf y}>}{\bf d}{\bf y}
$$
where $\L^*\in A^*_{\R}$ is the dual cone to $\L$ and 
${\bf d}{\bf y}$ is the Lebesgue measure on $A^*_{\R}$ 
normalized by the dual lattice $A^*\subset A^*_{\R}$ where 
$A^* := {\rm Hom}(A, {\bf Z})$.  
}
\end{dfn}

\begin{rem}{\rm 
If  $\L$ is a finitely generated
rational polyhedral cone the function ${\cal X}_{\L}({\bf s})$
is a rational function in ${\bf s}$. However, the explicit 
determination of this function might pose serious computational
problems. 
}
\end{rem}

\begin{dfn}
{\rm Let $V$ be an 
${\cal L}$-primitive smooth quasi-projective algebraic variety
with a metrized invertible sheaf ${\cal L}$ and $\a_{\cal L}(V)>0$.
Let  $X$ be any resolution of singularities 
$\rho\;:\; X\ra \overline{V}^{\cal L}$. 
We consider the triple
$$
({\rm Pic}(V,{\cal L}),{\rm Pic}(V,{\cal L})_{\R},
\L_{\rm eff}(V,{\cal L}))
$$
and the corresponding ${\cal X}$-function. 
Assuming that ${\rm Pic}(V, {\cal L})$ is a finitely generated 
abelian group (cf. \ref{vanish}) and that 
$\L_{\rm eff}(V,{\cal L})$ is a polyhedral 
cone (cf. \ref{polyhed}), we define the constant $\gamma_{\cal L}(V)\in \Q$ by
$$
\gamma_{\cal L}(V):= 
{\cal X}_{\L_{\rm eff}(V,{\cal L})}(\tilde{\rho}(-[K_X])).
$$
\label{gamma-dfn}
}
\end{dfn}

\subsection{${\cal L}$-primitive fibrations and descent of metrics}

Let $V$ be an ${\cal L}$-primitive variety of dimension $n$. 
We show that there exists a canonical measure 
on $V({\bf C})$ which is
uniquely defined up to a positive constant. 

In order to construct this measure we choose a 
resolution of singularities 
$\rho \; : \; X  \rightarrow   \overline{V}^{\cal L}$   and 
a positive integer $k_2$ such that 
$k_2D$ is a Cartier divisor, where  ${\cal O}(D) \cong 
(\rho^*L)^{\otimes \alpha_{\cal L}(V) } 
\otimes K_X$. Then $k_1 = k_2 \alpha_{\cal L}(V)$ is 
a positive integer. Let $g \in \H^0(X, {\cal O}(k_2D))$ be 
a non-zero global section (by \ref{rigid2}, 
it is uniquely defined up to 
a non-zero constant).  
We define a measure ${\bf \omega}_{{\cal L}}(g)$ 
on $X({\bf C})$ as follows. 
Choose local complex analytic coordinates 
$z_{1}, \ldots, z_{n}$ in 
some open neighborhood $U_x \subset X({\bf C})$ of a point 
$x \in V({\bf C})$.  
We write the restriction of the global section $g$  to $U_x$ as 
\[ g = s^{k_2\alpha_{\cal L}(V)} 
( dz_{1} \wedge  \cdots \wedge  dz_{n})^{\otimes k_2} = 
s^{k_1} ( dz_{1} \wedge  \cdots \wedge  dz_{n})^{\otimes k_2}, \]
where $s$ is a local section of $L$ . Then we set  
\[ {\bf \omega}_{{\cal L}}(g) : = \left(\frac{\sqrt{-1}}{2} \right)^n 
\|s\|^{ \alpha_{\cal L}(V)}  
(dz_1 \wedge d\overline{z}_1) \wedge  
   \cdots \wedge (dz_n \wedge d\overline{z}_n). \]
By a standard argument, one obtains that ${\bf \omega}_{{\cal L}}(g)$ 
doesn't depend on the choice 
of local coordinates in $U_x$ and  
that it extends to  the whole complex space 
$X({\bf C})$. It remains to notice that the restriction of 
the measure  ${\bf \omega}_{{\cal L}}(g)$ to $V({\bf C}) \subset X({\bf C})$ 
does not depend on the choice of $\rho$. So we obtain a well-defined 
measure on $V({\bf C})$. 

\begin{rem}
{\rm We note that the measure  $\omega_{\cal L}(g)$ depends on the 
choice of $g \in \H^0(X, {\cal O}(k_2D))$. More precisely, it multiplies 
by $|c|^{1/k_2}$ if we multiply $g$ by some non-zero complex number $c$. 
Thus we obtain that  
the mapping
\[ \| \cdot \|_{\cal L} \; : 
\; \H^0(X, {\cal O}(k_2D)) \rightarrow {\bf R}_{\geq 0} \]
\[ g \mapsto \int_{X({\bf C})} \omega_{\cal L}(g) = 
\int_{V({\bf C})} \omega_{\cal L}(g) \]
satisfies the property 
\[ \| c g \|_{\cal L} = |c|^{1/k_2}  \| g \|_{\cal L}\;\; \forall 
c \in {\bf C}^*. \]
}
\end{rem}

\begin{dfn}
{\rm Let $V$ be a smooth quasi-projective variety with 
an ample metrized invertible sheaf ${\cal L}$, 
 $\rho \, : \,X  \rightarrow   \overline{V}^{\cal L}$ a resolution 
of singularities. A regular 
projective morphism $\pi\, :\, X \ra Y$ to a projective 
variety $Y$ $({\rm dim}\, Y < {\rm dim}\, X)$ 
is called an ${\cal L}$-{\bf primitive fibration on} $V$  
if there exists a Zariski dense  open subset 
$U \subset Y$ such that the following conditions are satisfied:

{(i)} for any point $y \in U({\bf C})$ the fiber  
$V_y = \pi^{-1}({y}) \cap V$ is a smooth quasi-projective 
${\cal L}$-primitive subvariety;

(ii) $\alpha_{\cal L}(V) = \alpha_{\cal L}(V_y) > 0$ for all 
$y \in U({\bf C})$; 

(iii)  for any $k \in {\bf N}$ such that 
$k\alpha_{\cal L}(V) \in {\bf Z}$, 
\[ L_k := {\rm R}^0\pi_* \left(\rho^*(L)^{\otimes \alpha_L(V)} \otimes K_X 
\right)^{\otimes k} \]
is an ample invertible sheaf on $Y$. 
\label{prim-fb}}
\end{dfn}

\noindent
We propose the following 
version of the Fibration Conjecture of   
Fujita (see \cite{fujita0}):

\begin{conj}
{\sc (Existence of Fibrations)} 
Let $V$ be an arbitrary  smooth quasi-projective variety with 
an ample metrized invertible sheaf ${\cal L}$ and $\alpha_{\cal L}(V)> 0$. 
Then there exists a resolution of singularities 
 $\rho \, : \,X  \rightarrow   \overline{V}^{\cal L}$ 
such that $X$ admits an ${\cal L}$-primitive 
fibration  $\pi\, :\, X \ra Y$ on some dense Zariski open subset 
$V' \subset V$. 
\label{conj-fb}
\end{conj}

\begin{rem}
{\rm From the viewpoint of the Minimal Model Program, 
Conjecture  \ref{conj-fb} 
is equivalent to the statement  about the conjectured 
existence of ${\bf Q}$-Fano fibrations 
for algebraic varieties of negative Kodaira-dimension 
(cf. \cite{KMM}). The existence of 
an ${\cal L}$-primitive fibration is equivalent to 
the fact that 
the graded algebra
\[ {\rm R}(V, {\cal L}) = 
\oplus_{\nu \geq 0} \H^0(X, M^{k \nu}), \;\; 
M:= \rho^*(L)^{\otimes \alpha_L(V)} \otimes K_X  \]
is finitely generated (cf. 2.4 in \cite{BaMa}). One can 
define $Y$ as ${\rm Proj}\,{\rm  R}(V, {\cal L})$ and 
$X$ as a common resolution 
of singularities of  $\overline{V}^{\cal L}$ and of the 
indeterminacy locus and generic fiber of the natural rational map 
$\overline{V}^{\cal L} \ra Y$ (cf. \cite{H}). 
 }
\label{proj-fb}
\end{rem}
 
It is important to observe 
that a  metric $\| \cdot \|$ on $L$ induces natural 
metrics on all ample invertible sheaves 
$L_k$ on $Y$: 

\begin{dfn}
{\rm Let $L_k$ be an ample invertible 
sheaf on $Y$ as above. 
We define a metric $\| \cdot \|_{{\cal L}, k}$  
on $L_k$ as follows.  
Let $y \in Y({\bf C})$ be a closed point, 
$U \subset Y$ be a Zariski  
open subset containing $y$, and 
$s \in \H^0(U, L_k)$ is a section with 
$s(y) \neq 0$. Then we set 
\[ \|s(y)\|_{{\cal L},k}  := 
\left( \int_{V_y({\bf C})} \omega_{\cal L}(\pi^* s) \right)^k , \]
where $V_y({\bf C})$ is the fiber over $y$ of the ${\cal L}$-primitive 
fibration $\pi^*$, $\pi^* s$ the $\pi$-pullback of $s$ restricted to 
${\cal L}$-primitive variety $V_y({\bf C})$, and 
 $\omega_{\cal L}(\pi^* s)$ the corresponding to  $\pi^* s$ volume 
measure on   $V_y({\bf C})$. 
We call  $\| \cdot \|_{{\cal L}, k}$ a $k$-{\bf adjoint descent} to $Y$ 
of a metric $\| \cdot \|$ on $L$. }
\end{dfn}

\section{Heights and asymptotic formulas} 

\subsection{Basic terminology and notations}

Let $F$ be a number field, ${\cal O}_F \subset F$ the 
ring of integers in $F$, ${\rm Val}(F)$ the set of all valuations of $F$, 
$F_v$ the completion of $F$ with respect to a valuation $v \in {\rm Val}(F)$, 
${\rm Val}(F)_{\infty} = \{ v_1, \ldots, v_r \}$ the set of all archimedean 
valuations of $F$. 
For any algebraic variety $X$ over 
a field $F$ we denote by 
$X(F)$ the  set of its $K$-rational points.

\begin{dfn}
{\rm Let $E$ be a vector space of dimension $m+1$ over $F$, 
${\cal O}_E \subset E$ a projective ${\cal O}_F$-module of rank $m+1$ 
and $\| \cdot \|_{v_1}, \ldots, \| \cdot \|_{v_r}$ the set 
of Banach norms on  the real or complex vector spaces 
$E_{v_i} = E \times_F F_{v_i}$ corresponding to elements of  
${\rm Val}(F)_{\infty} = \{ v_1, \ldots, v_r \}$. It is well-known that  
the above data for $E$ define a 
family $\{ \| \cdot \|_v, \; v \in {\rm Val}(F) \}$ of $v$-adic metrics 
for  a standard invertible sheaf  ${\cal O}(1)$ on  ${\bf P}(E)$. 
If   $x \in X(F)$ is a point and  
$s \in \H^0(U, {\cal O}(1))$ is a section  over an open subset  
$U  \subset {\bf P}(E)$ 
containing $x$, then  we denote 
by  $\|s(x)\|_v$ the corresponding 
$v$-adic norm of 
$s$ at $x$. We set 
$\tilde{\cal O}(1) = 
( {\cal O}(1), \|\cdot \|_v) $ to be the  
standard invertible sheaf ${\cal O}(1)$ on ${\bf P}(E)$ 
together with a family $v$-adic metrics  
$\|\cdot \|_v$ defined 
by the above data and we
call $\tilde{\cal O}(1)$ the 
{\bf standard ample metrized invertible 
sheaf} on ${\bf P}(E)$. }
\end{dfn} 

\begin{dfn}
{\rm  Let $X$ be an  algebraic variety over $F$. For any point   
 $x \in X(F)$ and any regular function   
$f \in \H^0(U, {\cal O}_X)$ on  an open subset  
$U  \subset {\bf P}(E)$ 
containing $x$, we define the $v$-adic norm  
$\|f(x)\|_v := |f(x)|_v$. We call 
this family of $v$-adic metrics on ${\cal O}_X$ 
the {\bf canonical metrization} of the structure sheaf ${\cal O}_X$. }
\end{dfn}

\begin{dfn}
{\rm Let $X$ be a quasi-projective algebraic variety,  
$L$ a very ample  invertible sheaf on $X$, 
$i \; :\; X \hookrightarrow  {\bf P}(E)$  an embedding with  
$L = i^* {\cal O}(1)$. We denote by ${\cal L} = (L, \|\cdot \|_v )$ 
the sheaf $L$ together with 
$v$-adic metrics  induced from a 
family of $v$-adic metrics on 
the standard ample metrized 
invertible sheaf  $\tilde{\cal O}(1)$ on 
${\bf P}(E)$. In this situation we call ${\cal L}$ 
 {\bf a very ample metrized sheaf} on $X$ and write 
 ${\cal L} = i^*\tilde{O}(1)$. 
}
\end{dfn}

\begin{dfn}
{\rm Let ${\cal L} =( L,  \|\cdot \|_v) $ be 
a very ample metrized invertible sheaf on $X$. Then 
for any point  $x \in X(F)$, the ${\cal L}$-{\bf height} of 
$x$ is defined as 
\[ H_{\cal L}(x) = \prod_{v \in {\rm Val}(F)} \|s(x)\|_v^{-1}, \]
where $s \in \Gamma(U, L)$ is a nonvanishing at $x$ 
section of $L$ over some open 
subset $U \subset X$. }
\end{dfn}

\begin{rem}
{\rm Using a canonical metrization of the structure sheaf ${\cal O}_X$,   
the  linear mapping 
\[ S^k(\Gamma(U, L)) \rightarrow  \Gamma(U, L^{\otimes k}) \;\; (k > 0), \]
and the $F$-bilinear mapping 
\[  \Gamma(U, L^{\otimes k}) \times \Gamma(U, L^{\otimes -k}) \rightarrow 
 \Gamma(U, {\cal O}_X), \] 
one immediately sees that a family of $v$-adic metrics 
on an invertible sheaf $L$ allows to define a family of $v$-adic 
metrics on $L^{\otimes k}$ and on 
any  invertible sheaf $M$ such that there exist 
integers $k_1, k_2$ $(k_2 \neq 0)$ 
with $L^{\otimes k_1} = M^{\otimes k_2}$. 
In this situation we write 
$$
{\cal M} =  (M,  \|\cdot \|_v) := (L^{\otimes{k_1/k_2}},  
\|\cdot \|_v^{k_1/k_2}),
$$ 
or simply $ {\cal M} ={\cal L}^{k_1/k_2}$. 
Obviously, 
one obtains 
\[ H_{\cal M}(x)  =  (H_{\cal L}(x))^{k_1/k_2}\; \; \;\forall \; 
x \in X(F). \]}
\label{fr-norms}
\end{rem}

\begin{dfn}
{\rm Let $L$ be an ample invertible sheaf on 
a quasi-projective variety $X$ 
and $k$ a positive integer such that $L^{\otimes k}$
is very ample.
We  define an 
{\bf ample metrized invertible sheaf}  ${\cal L}=( L,  \|\cdot \|_v)$ 
on $X$ associated with $L$ by considering 
${\cal L}^{\otimes k}: = (L^{\otimes k}, \| \cdot \|_v^k)$ as  
a very ample metrized invertible sheaf on $X$.
}
\end{dfn}

\subsection{Weakly and strongly ${\cal L}$-saturated varieties}

Let $V$ be an arbitrary quasi-projective algebraic variety over $F$ with 
an ample metrized invertible sheaf ${\cal L}$. We always assume 
that $V(F)$  is infinite and  set 
$$
N(V,{\cal L}, B):=\#\{x\in V(F)\; : \; H_{\cal L}(x)\le B\}.
$$
Here and in 3.4 we will work under the following 

\begin{assume}
{\rm For all quasi-projective $V',V$  with $V'\subset V$
and $|V(F)|=\infty$ there exists the limit
$$
\lim_{B\ra \infty} \frac{N(V',{\cal L}, B)}{N(V,{\cal L}, B)}.
$$
}
\label{assumption}
\end{assume}

\begin{dfn}
{\rm We call an irreducible  
quasi-projective algebraic variety $V$ 
with an ample metri\-zed 
invertible sheaf ${\cal L}$ 
{\bf weakly  ${\cal L}$-saturated} if 
for any Zariski locally closed subset $W \subset V$ with 
${\rm dim}\, W < {\rm dim}\, V$, one has 
\[ {\lim}_{B \rightarrow \infty} 
\frac{N(W,{\cal L},B)}{N(V,{\cal L},B)} < 1. \]
}
\end{dfn}

\begin{dfn}
{\rm We call an irreducible 
quasi-projective algebraic variety $V$ 
with an ample  metri\-zed 
invertible sheaf ${\cal L}$ 
{\bf strongly  ${\cal L}$-saturated} if 
for any dense Zariski open subset $U \subset V$, one 
has 
\[ {\lim}_{B \rightarrow \infty} 
\frac{N(U,{\cal L},B)}{N(V,{\cal L},B)} = 1. \]
}
\end{dfn}

\begin{dfn}
{\rm Let $V$ be a weakly ${\cal L}$-saturated 
variety, $W \subset V$ a locally 
closed strongly saturated 
subvariety of smaller dimension.  
Then we call $W$ an  ${\cal L}$-{\bf target} of $V$, if 
$$
0 < {\lim}_{B \rightarrow \infty} 
\frac{N(W,{\cal L},B)}{N(V,{\cal L},B)} < 1. 
$$
}
\end{dfn}

\begin{theo}
Let $V$ 
be an arbitrary  quasi-projective 
algebraic variety with an 
ample metrized invertible sheaf ${\cal L}$. Assume that 
$|V(F)|= \infty$ and that \ref{assumption} holds. Then
we have:

{\rm (i)} if $V$ is  strongly   
${\cal L}$-saturated then 
$V$ is weakly  ${\cal L}$-saturated; 

{\rm (ii)} $V$ contains finitely many 
weakly ${\cal L}$-saturated subvarieties  
$W_1, \ldots, W_k$ with 
\[ {\lim}_{B \rightarrow \infty} 
\frac{N(W_1 \cup \cdots 
\cup W_k,{\cal L}, B)}{N(V,{\cal L},B)} =  1. \]

{\rm (iii)} $V$ 
contains a strongly  ${\cal L}$-saturated 
subvariety $W$ having the property 
\[ 0 < {\lim}_{B \rightarrow \infty} 
\frac{N(W,{\cal L},B)}{N(V,{\cal L},B)} < 1. \]

{\rm (iv)} if $V$ is weakly saturated and if 
it doesn't contain a dense Zariski open
subset $U \subset V$ which is 
strongly saturated then $V$ contains infinitely many 
${\cal L}$-targets. 
\label{ws-sat}
\end{theo}

\noindent
{\em Proof.} (i) Let $W \subset V$ be  a  Zariski closed 
subset with ${\rm dim}\, W < {\rm dim}\, V$ and $U = V \setminus W$. Then 
\[ {\lim}_{B \rightarrow \infty} 
\frac{N(W,{\cal L},B)}{N(V,{\cal L},B)} + 
{\lim}_{B \rightarrow \infty} 
\frac{N(U,{\cal L},B)}{N(V,{\cal L},B)} = 
{\lim}_{B \rightarrow \infty} 
\frac{N(V,{\cal L},B)}{N(V,{\cal L},B)} =1.
\]
Since $V$ is strongly ${\cal L}$-saturated, we have 
\[ {\lim}_{B \rightarrow \infty} 
\frac{N(U,{\cal L},B)}{N(V,{\cal L},B)} =1 \]
and therefore
\[ {\lim}_{B \rightarrow \infty} 
\frac{N(W,{\cal L},B)}{N(V,{\cal L},B)} =0 <1. \]

(ii) Let $W \subset V$ be a  
minimal  Zariski closed subset  such that 
\[ {\lim}_{B \rightarrow \infty} 
\frac{N(W,{\cal L},B)}{N(V,{\cal L},B)} =  1 \]
and $W_1, \ldots, W_k$ irreducible components of $W$. 
It immediately follows from the minimality of $W$ that 
each $W_i$ is weakly  saturated.

(iii) Let $W \subset V$ be an 
irreducible Zariski closed subset  of 
minimal dimension  
such that 
\[ {\lim}_{B \rightarrow \infty} 
\frac{N(W,{\cal L},B)}{N(V,{\cal L},B)} =1. \]
The minimality of $W$ implies that 
$W$ is weakly  saturated. 

(iv) By (iii) the set of ${\cal L}$-targets is nonempty. 
Assume that the set of all ${\cal L}$-targets is finite: 
$\{W_1, \ldots, W_k \}$.  
The strong saturatedness of each $W_i$ implies that 
\[ {\lim}_{B \rightarrow \infty} \frac{N(W_{i_1} \cap \cdots 
\cap W_{i_l},{\cal L}, B)}{N(V,{\cal L}, B)} = 0 \]
for all pairwise different $i_1, \ldots, i_l \in \{1,\ldots, k\}$ 
and $l \geq 2$. In particular, 
one has 
\[ \sum_{i=1}^k  {\lim}_{B \rightarrow \infty} 
\frac{N(W_i,{\cal L},B)}{N(V,{\cal L}, B)} =  
{\lim}_{B \rightarrow \infty} 
\frac{N(W_1 \cup \cdots \cup W_k,{\cal L},B)}{N(V,{\cal L},B)} < 1. \]
We set $U:= V \setminus (W_1 \cup \cdots \cup W_k)$. Then 
\[  {\lim}_{B \rightarrow \infty} 
\frac{N(U,{\cal L},B)}{N(V,{\cal L},B)} > 0. \]
Since $U$ is not strongly saturated, 
there exists an irreducible Zariski closed subset $W_0 \subset V$ 
of minimal dimension $< {\rm dim}\, V$ such that 
\[ {\lim}_{B \rightarrow \infty} 
\frac{N(W_0\cap U,{\cal L},B)}{N(U,{\cal L},B)} > 0. \]
It follows from the minimality of $W_0$ that $W_0$ is strongly 
saturated.
 On the other hand,  one has 
\[  {\lim}_{B \rightarrow \infty} 
\frac{N(W_0,{\cal L},B)}{N(V,{\cal L},B)} \geq 
 {\lim}_{B \rightarrow \infty} 
\frac{N(W_0\cap U,{\cal L},B)}{N(V,{\cal L},B)} > 0, \]
i.e., $W_0$ is an ${\cal L}$-target and $W_0 \not\in 
\{ W_1, \ldots, W_k \}$. Contradiction.  
\hfill $\Box$

\begin{dfn}
{\rm Let $V$ be a  weakly ${\cal L}$-saturated variety 
and $W_1, W_2, \ldots$ an infinite sequence 
of strongly saturated irreducible subvarieties 
$W_i$ having the property 
\[ 0 < \theta_i := 
{\lim}_{B \rightarrow \infty} 
\frac{N(W_i,{\cal L},B)}{N(V,{\cal L},B)} < 1\;\; \forall i > 0. \]
We say that the set $\{W_1, W_2, \ldots  \}$ 
forms an {\bf   
asymptotic arithmetic ${\cal L}$-fibration} on $V$,  
if the following equality holds
\[ \sum_{i=1}^{\infty} \theta_i = 1. \] } 
\end{dfn} 
\medskip

We expect that the main source of 
examples of weakly and strongly 
saturated varieties should come 
from the following situation:

\begin{prop}
{\rm  Assume \ref{assumption} 
and let $G \subset PGL(n+1)$ 
be a connected linear 
algebraic group acting on $\P^n$ and  
$V :=Gx \subset \P^n$ a 
$G$-orbit of a point $x \in \P^n(F)$. 
Then  $V$ is weakly $\tilde{O}(1)$-saturated.
\label{sat1}
}
\end{prop}

\noindent 
{\em Proof.} Let $W \subset V$ be an arbitrary 
locally closed 
subset with ${\rm dim}\, W < {\rm dim}\, V$, 
$\overline{W} \subset V$ 
its Zariski closure in $V$ and 
$ U := V \setminus \overline{W} \subset V$ 
the corresponding 
dense  Zariski open 
subset of $V$. Then $V$ is covered by the 
open subsets $gU$, 
where $g$ runs over all elements in $G(F)$ 
(this follows from the fact
that $G$ is unirational
and that $G(F)$ is Zariski dense in $G$ \cite{borel}). 
Therefore, the orbit of $x\in V(F)$ 
under $G(F)$ is Zariski dense in $V$.
Since the Zariski topology is noetherian we 
can choose a finite subcovering:
$V = \bigcup_{i =1 }^k g_iU$ ($g_i$ in $G(F)$). 
Considering  $g_i\in G(F)$ as matrices 
in $PGL(n+1)$ and using standard 
properties of heights \cite{lang}, 
one obtains positive constants
$c_i$ such that 
$$
 H_{\cal L}(g_i(x)) \le c_iH_{\cal L}(x)
$$
for all $x\in \P^n(F)$. 
It is clear that
for $c_0: =\sum_{i =1}^k c_i$ we have
$$
N(U,{\cal L},B) \le N(V, {\cal L},B) \le c_0 N(U, {\cal L},B).
$$
It follows that
$$
\frac{1}{c_0} \leq 
{\lim}_{B \rightarrow \infty} 
\frac{N(U, {\cal L},B)}{N(V,{\cal L},B)} \leq 1.
$$
Hence 
$$
{\lim}_{B \rightarrow \infty} 
\frac{N(W, {\cal L},B)}{N(V,{\cal L},B)} \leq  
{\lim}_{B \rightarrow \infty} 
\frac{N(\overline{W}, {\cal L},B)}{N(V,{\cal L},B)}
\leq  1 -\frac{1}{c_0} <1.
$$

\vskip 0,5cm
\noindent
We can reformulate the statement of \ref{sat1} as follows: 

\begin{prop} {\rm 
Let $G$ be a connected linear algebraic group, $H \subset G$ a
closed subgroup and  $V : = G/H$.  If
\ref{assumption} holds then $V$ is weakly saturated 
with respect to
any $G$-equivariant projective embedding of $V$.
\label{equiv-sat}
}
\end{prop}

\noindent
It is easy to  see   that $V = G/H$ is not necessarily weakly saturated 
with respect to projective embeddings which are not
$G$-equivariant:
 
\begin{exam}
{\rm Let $S \subset \P^8$ be the anticanonically embedded 
Del Pezzo surface  which is a blow up of a
rational point in $\P^2$. Denote by ${\cal L}$ the metrized anticanonical 
sheaf on $S$. 
The unique exceptional curve $C \subset S$ is
contained in the union of two open subsets $U_0, U_1 \subset S$
where $U_0 \cong U_1 \cong {\bf A}^2$. Therefore, $S$  can be
considered as a projective compactification of the algebraic
group ${\bf G}_a^2$ (after an identification of ${\bf G}_a^2$ with 
 $U_0$ or $U_1$). This compactification is not ${\bf G}_a^2$-equivariant. 
One has  
$$
a_{\cal L}(U_0) = a_{\cal L}(U_1) = a_{\cal L}(C) = 2, 
$$
but 
$$
a_{\cal L}(U_0 \setminus C) = a_{\cal L}(U_1 \setminus C) =1.
$$
Hence, $U_0$ and  $U_1$ are not weakly 
${\cal L}$-saturated.
 }
\end{exam}

It is easy to show that an equivariant compactification of $G/H$ 
is not necessarily strongly ${\cal L}$-saturated: 

\begin{exam}
{\rm Let $V = {\bf P}^1 \times {\bf P}^1$. Then $V$ is a $G$-homogeneous 
variety with $G = GL(2) \times GL(2)$. However, $V$ is not strongly 
${\cal L}$-saturated for $L := \pi_1^*{\cal O}(k_1) \otimes   
\pi_2^*{\cal O}(k_2)$ ($k_1,k_2  \in {\bf N}$), if $k_1 \neq k_2$. } 
\end{exam}

\subsection{Adelic ${\cal L}$-measure and $\tau_{\cal L}(V)$}

Now we define an adelic 
measure ${\bf \omega}_{\cal L}$ corresponding to an 
ample metrized  invertible sheaf ${\cal L}$ on an ${\cal L}$-primitive 
variety $V$ with $\alpha_{\cal L}(V) >0$ which satisfies
the assumption \ref{vanish}. This is a generalization of 
a construction due to Peyre (\cite{peyre}) for $V$ 
being a smooth projective variety and ${\cal L}$ the metrized canonical
line bundle, which in its turn is a generalization of the classical 
construction of Tamagawa measures on  the adelic points of algebraic groups.

Let  $V$ be an ${\cal L}$-primitive variety of dimension $n$, 
$\rho \, : \, X  \rightarrow   \overline{V}^{\cal L}$ a resolution of 
singularities, $k_2$ a positive integer  such that $k_1 = 
k_2 \alpha_{\cal L}(V)  \in {\bf Z}$ and 
$$
( \rho^*(L)^{\alpha_{\cal L}(V)} \otimes K_X)^{\otimes k_2} \cong  
{\cal O}(D),
$$
where $D$  
is a rigid effective Cartier divisor 
on $X$.    

\begin{dfn}
{\rm Let $g$ be a non-zero element of the 
$1$-dimensional $F$-vector space $\H^0(X, {\cal O} (D))$  
($g$  is defined uniquely up to an element of $F^*$).  
Let $v  \in {\rm Val}(F)$. We define a measure ${\bf \omega}_{{\cal L},v}(g)$ 
on $V(F_v)$ as follows. 
Choose local $v$-analytic coordinates 
$x_{1,v}, \ldots, x_{n,v}$ in 
some open neighborhood $U_x \subset X(F_v)$ of a point 
$x \in X(F_v)$.  
We write the restriction of the global section $g$  to $U_x$ as 
\[ g = s^{k_1} ( dx_{1,v} \wedge  \cdots \wedge  dx_{n,v})^{k_2} \]
where $s$ is a local section of $L$. Define a $v$-adic measure on $U_x$
as 
\[ {\bf \omega}_{{\cal L},v}(g) : = \|s\|_v^{k_1/k_2}  dx_{1,v}
   \cdots  dx_{n,v} =  \|s\|_v^{ \alpha_{\cal L}(V)}  dx_{1,v}
   \cdots dx_{n,v}, \]
where $dx_{1,v}\cdots dx_{n,v}$ is the usual normalized
Haar measure on $F_v^n$.
By a standard argument, one obtains that 
${\bf \omega}_{{\cal L},v}(g)$ 
doesn't depend on the choice 
of local coordinates in $U_x$ and  
that it extends to  the whole $v$-adic space 
$X(F_v)$. The restriction of  
${\bf \omega}_{{\cal L},v}(g)$ doesn't depend 
on the choice of $\rho$. So we obtain a 
well-defined $v$-adic measure 
on $V(F_v)$. }
\end{dfn}

\begin{rem}
{\rm We remark that ${\bf \omega}_{{\cal L},v}(g)$ 
depends on the choice of a  global section 
$g \in \H^0(X, {\cal O} (D))$: 
if $g' = cg$ $( c\in F^*)$ is another 
global section, then 
\[  {\bf \omega}_{{\cal L},v}(g') =|c|_v^{1/k}
{\bf \omega}_{{\cal L},v}(g). \]
}
\end{rem}

\noindent
Our next goal is to obtain an 
explicit formula for the integral 
\[ d_v(V):  = \int_{V(F_v)}{\bf \omega}_{{\cal L},v}(g) = 
 \int_{X(F_v)}{\bf \omega}_{{\cal L},v}(g)   \]
for almost all $v \in {\rm Val}(F)$. We use the $p$-adic integral formula 
of Denef (\cite{denef} Th. 3.1) in the same way  as the $p$-adic formula 
of A. Weil (\cite{weil} Th. 2.2.3) was used by Peyre in \cite{peyre}.

\begin{rem}
{\rm By \cite{AW,H}, we can choose $\rho$ in such a way that 
$\rho$  is defined over $F$ and all irreducible components $D_1, \ldots, 
D_l$ $(l \geq 0)$  of the support 
of $D = \sum_{i =1}^l m_i D_i$  are smooth divisors 
with normal crossings over the algebraic closure 
$\overline{F}$.}
\end{rem}

\begin{dfn}
{\rm  Let $G$ be the image of ${\rm Gal}(\overline{F}/F)$ 
in the symmetric group $S_l$ that acts by permutations on  $D_1, \ldots, 
D_l$. We set $I: = \{ 1, \ldots, l\}$ and 
denote by $I/G$ the set of 
all $G$-orbits in $I$. 
For any $J  \subset I$ we set   
\[ D_J : = \left\{ \begin{array}{ll} 
\bigcap_{j \in J} D_j, \; &  \;  \mbox{if $J \neq \emptyset$} \\
 X \setminus \bigcup_{j \in I} D_j, \; &  \; \mbox{if $J =  \emptyset$} 
 \end{array} \right. ,  \]
\[ D_J^{\circ} = D_J \setminus \bigcup_{j \not\in J} D_j.  \]
($D_J$ is defined over $F$ if and only if $J$ is 
a union of some of $G$-orbits in  $I/G$.). }  
\end{dfn}

We can extend $X$ to a projective 
scheme ${\cal X}$ of finite type over ${\cal O}_F$ and divisors 
$D_1, \ldots, D_l$ to codimension-$1$ subschemes ${\cal D}_1, 
\ldots, {\cal D}_l$ in ${\cal X}$ such that 
for almost all non-archimedean $v \in {\rm Val}(F)$ 
the reductions of $X$  and ${\cal D}_1, \ldots, {\cal D}_l$ 
modulo $\wp_v \subset {\cal O}_F$ are smooth projective varieties 
${\cal X}_v$ and ${\cal D}_{v,1}, \ldots, {\cal D}_{v,l}$ over 
the algebraic closure $\overline{k_v}$ of the 
residue field $k_v$ with  ${\cal D}_{v,i} \neq  {\cal D}_{v,j}$ 
for $i \neq j$.  Moreover, we can assume that 
$$
{\cal D}_{v,J} : = \left\{ \begin{array}{ll} 
\bigcap_{j \in J} {\cal D}_{v,j}, \; &  \;  \mbox{if $J \neq \emptyset$} \\
 X \setminus \bigcup_{j \in I} {\cal D}_{v, j}, \; &  \;
\mbox{if $J =  \emptyset$} 
 \end{array} \right.
$$
are also smooth over  $\overline{k_v}$. 

\begin{dfn}
{\rm A non-archimedean valuation $v \in {\rm Val}(F)$
which satisfies all the
above assumptions will be called 
a {\bf good valuation} for the pair $({\cal X}, 
\{ {\cal D}_i \}_{i \in I})$. }
\end{dfn} 

\begin{dfn}
{\rm Let $G_v \subset G$ be a cyclic subgroup generated by a representative 
of the Frobenius element in  ${\rm Gal}(\overline{k_v}/k_v)$. We denote 
by $I/G_v$ the set of all $G_v$-orbits in $I$. If $j \in I/G_v$, then 
we set  $b_j$ to be the length of the corresponding $G_v$-orbit 
and put  $r_j = m_j/k_2$, where $m_j$ is the multiplicity of irreducible 
components of $D$ corresponding to the $G_v$-orbit $j$. }
\end{dfn}

The following theorem is a slightly generalized version of 
Th. 3.1 in \cite{denef}:  

\begin{theo}
Let $v \in {\rm Val}(F)$ be a good non-archimedean valuation 
for  $({\cal X}, \{ {\cal D}_i \}_{i \in I})$. Then 
\[ d_v(V) = \int_{X(F_v)} {\bf \omega}_{{\cal L},v}(g) = 
\frac{c_{\emptyset}}{q_v^n} +  \frac{1}{q_v^n} 
 \sum_{\emptyset \neq J \subset I/{G_v} } c_J \prod_{j \in J} 
\left( \frac{q_v^{b_j} -1}{q_v^{b_j(r_j + 1)} -1} \right), \]
where $q_v$ is the cardinality of $k_v$, and 
$c_J$ is the cardinality  of the set of $k_v$-rational points
in ${\cal D}^{\circ}_J$. 
\end{theo}

Let us consider an exact sequence of 
Galois ${\rm Gal}(\overline{F}/F)$-modules: 
\begin{equation}
 0 \ra {\bf Z}[{\cal D}_1] \oplus \cdots \oplus  
{\bf Z}[{\cal D}_l] \ra {\rm Pic}({\cal X}) 
\stackrel{\tilde{\rho}}{\ra}  {\rm Pic}({\cal X} \setminus 
\bigcup_{i =1}^l {\cal D}_i ) {\ra} 0 
\label{sh-3}
\end{equation}

\begin{theo}
Assume that $X$ has  the property  ${\rm h}^1(X, {\cal O}_X) =0$. 
Then 
\[ d_v(V) = 1 + \frac{1}{q_v} {\rm Tr}(\Phi_v | {\rm Pic}({\cal X} \setminus 
\bigcup_{i =1}^l {\cal D}_i )\otimes {\bf Q}_l) + 
O\left(\frac{1}{q_v^{1+ \varepsilon}} \right), \]
where
\[ \varepsilon = \min \{ 1/2, r_1, \ldots, r_l \} \]
and $\Phi_v$ is the Frobenius morphism. 
\end{theo}

\noindent
{\em Proof.} 
By conjectures of Weil proved 
by Deligne \cite{deligne}, one has 
\[ \frac{c_J}{q_v^n} = O\left(\frac{1}{q_v} \right) \; J \neq \emptyset \]
(since ${\rm dim}\, D_J \leq n-1$ for $J \neq \emptyset$)
and 
\[  \frac{c_{\emptyset}}{q_v^n}  
= \sum_{k =0}^{2n} (-1)^k {\rm Tr}(\Phi_v | \H^k_c({\cal X} \setminus 
\bigcup_{i =1}^l {\cal D}_i, {\bf Q}_l)), \]
where  $\H^k_c(\cdot , {\bf Q}_l)$ denotes the \'etale cohomology group
with compact supports. 
Using long cohomology sequence of the pair $({\cal X}, \bigcup_i {\cal D}_i)$, 
one obtains  isomorphisms
\[  \H^{k}_c({\cal X} \setminus 
\bigcup_{i =1}^l {\cal D}_i, {\bf Q}_l) =  \H^{k}_c({\cal X}, {\bf Q}_l)\;\; 
\mbox{\rm for $k = 2n, 2n-1$} \]
and the short exact sequence
\[ 0 \ra \H^{2n-2}_c({\cal X} \setminus 
\bigcup_{i =1}^l {\cal D}_i, {\bf Q}_l) \ra 
\H^{2n-2}_c({\cal X}, {\bf Q}_l) \ra \bigoplus_{i =1}^l 
\H^{2n-2}_c({\cal D}_i, {\bf Q}_l) \ra 0. \]
Using isomorphisms 
\[ \H^{2n-2}_c({\cal D}_i, {\bf Q}_l(n-1)) \cong {\bf Q}_l,  
\H^{2n}_c({\cal X}, {\bf Q}_l(n)) \cong {\bf Q}_l, \]  
Poincar\'e duality 
\[ \H^{2n-2}_c({\cal X}, {\bf Q}_l(n-1)) 
\times {\rm Pic}({\cal X}) 
\stackrel{\sim}{\ra} {\bf Q}_l \]
and the vanishing property  
${\rm h}^1(X, {\cal O}_X) =0$, we obtain 
\[   \H^{2n-1}_c({\cal X} \setminus 
\bigcup_{i =1}^l {\cal D}_i, {\bf Q}_l) = 0 \]
and 
\[ \frac{c_{\emptyset}}{q_v^n} = 1 + 
\frac{1}{q_v} {\rm Tr}(\Phi_v | 
{\rm Pic}({\cal X} \setminus 
\bigcup_{i =1}^l {\cal D}_i )\otimes {\bf Q}_l) 
 \]
\[ + \sum_{k =0}^{2n-3} \frac{ (-1)^k}{q_v^n}  
{\rm Tr}(\Phi_v | \H^k_c({\cal X} \setminus 
\bigcup_{i =1}^l {\cal D}_i, {\bf Q}_l)) \]
\[ = 1 + \frac{1}{q_v} {\rm Tr}(\Phi_v | 
{\rm Pic}({\cal X} \setminus 
\bigcup_{i =1}^l {\cal D}_i )\otimes {\bf Q}_l) + 
O\left(\frac{1}{q_v^{3/2}} \right). \]
On the other hand, for $J \neq \emptyset$ one has 
\[  \prod_{j \in J} 
\left( \frac{q_v^{b_j} -1}{q_v^{b_j(r_j + 1)} -1} \right) 
= O\left(\frac{1}{q_v^{ 1 + \varepsilon_0}} \right), \]
where $\varepsilon_0 = \min \{ r_1, \ldots, r_l \}$. 

\hfill $\Box$ 

\begin{dfn}
{\rm We define the convergency factors
\[ \lambda_v : = 
\left\{ \begin{array}{ll} L_v(1, {\rm Pic}({\cal X} \setminus 
\bigcup_{i =1}^l {\cal D}_i )), \;& \;  \mbox{\rm if $v$ is good} \\
1    \; & \; \mbox{\rm otherwise} \end{array} \right. \]
where $L_v$ is the local factor of the Artin $L$-function
corresponding to the $G$-module ${\rm Pic}({\cal X} \setminus 
\bigcup_{i =1}^l {\cal D}_i)$ 
and we set 
\[ \omega_{{\cal L},S} : = \sqrt{|{\rm disc}(F)|}^{-n} \prod_{v \in {\rm Val}(F)} 
\lambda_v^{-1} {\bf \omega}_{{\cal L},v}(g), \]
where ${\rm disc}(F)$ is the
absolute discriminant of ${\cal O}_F$
and $S$ is the set of bad valuations. 
}
\end{dfn}

By the product formula, 
$\omega_{{\cal L},S}$ doesn't 
depend on the choice of 
$g$. 

\begin{dfn}
{\rm
Denote by ${\bf A}_F$ the adele ring of $F$.
Let $\overline{X(F)}$ be the closure of $X(F)$ in 
$X({\bf A}_F)$ (in direct product topology).  
Under the vanishing assumption ${\rm h}^1(X,{\cal O}_X) =0$, 
we define the constant 
\[ \tau_{\cal L}(V) = 
\lim_{s \ra 1}(s-1)^{{\beta}_{\cal L}(V)} L_S(s, {\rm Pic}({\cal X} \setminus 
\bigcup_{i =1}^l {\cal D}_i )) \int_{\overline{X(F)}} \omega_{{\cal L},S}, \]
where  ${\beta}_{\cal L}(V)$ is  the rank of the submodule 
of ${\rm Gal}(\overline{F}/F)$-invariants of the module
${\rm Pic}({\cal X} \setminus 
\bigcup_{i =1}^l {\cal D}_i )$. 
\label{tau-dfn}
}
\end{dfn}

\subsection{Main strategy} 

Now we proceed to  discuss our 
main strategy in understanding the asymptotic 
for the number $N(V, {\cal L}, B)$ as 
$B \to \infty$ for an arbitrary 
${\cal L}$-polarized quasi-projective variety. 
Again, we shall make the assumption \ref{assumption}. 
Our approach consists 
in $4$  steps including  $3$ subsequent simplifications 
of the situation: 
\medskip

\noindent
{\bf Step 1 (reduction 
to weakly ${\cal L}$-saturated varieties):}  
By \ref{ws-sat} (ii), 
every quasi-projective ${\cal L}$-polarized  
variety $V$ contains a finite number 
of weakly ${\cal L}$-saturated 
varieties $W_1, \ldots W_k$ such that 
\[ 
{\lim}_{B \rightarrow \infty} 
\frac{N(W,{\cal L},B)}{N(V,{\cal L},B)} =  1.
\]

Therefore, it would be enough to understand 
separately the asymptotics of 
$N(W_i, {\cal L}, B)$, $i \in \{1, \ldots, k\}$
 modulo the asymptotics 
$N(W_{i_1} \cap \cdots \cap W_{i_l}, {\cal L}, B)$ 
for low-dimensional 
subvarieties $W_{i_1} \cap \cdots \cap W_{i_l}$, where 
$\{i_1, \ldots, i_l\} \subset \{1, \ldots, k\}$ 
are subsets of pairwise different 
elements with $l \geq 2$.  
\medskip

For  our next reduction step, we need: 

\begin{conj} 
Let $V$ be a weakly ${\cal L}$-saturated  
variety which doesn't contain an open Zariski 
dense and 
strongly  ${\cal L}$-saturated subset $U \subset V$. 
Then the 
set of ${\cal L}$-targets of $V$ forms 
an asymptotic arithmetic 
${\cal L}$-fibration.
\label{aaf} 
\end{conj}

\noindent
{\bf Step 2 (reduction to strongly ${\cal L}$-saturated varieties):} 
Let $V$ be an arbitrary  weakly ${\cal L}$-saturated  variety. Then 
either $V$ contains a strongly ${\cal L}$-saturated Zariski open subset 
or, according to \ref{aaf},  
we obtain   an asymptotic arithmetic 
${\cal L}$-fibration of $V$ by ${\cal L}$-targets. 
In the first situation, 
it is 
enough to understand the asymptotic of 
$N(U, {\cal L}, B)$ for the 
strongly ${\cal L}$-saturated variety $U$ 
(we note that the  complement $V \setminus U$ 
consists of low-dimensional irreducible components). 
In the second situation, 
it is 
enough to understand the asymptotic of 
$N(W_i, {\cal L}, B)$ for each of the
${\cal L}$-targets $W_i \subset V$. 
\medskip

For  our next reduction step, we need:

\begin{conj} 
Let $V$ be a smooth strongly 
${\cal L}$-saturated quasi-projective 
variety. Then the complex analytic 
variety $V({\bf C})$ is ${\cal L}$-primitive.
\end{conj} 

\noindent
{\bf Step 3 
(reduction to ${\cal L}$-primitive varieties):} 
Every quasi-projective algebraic 
variety $V$ is a disjoint union of finitely many 
locally closed smooth subvarieties $V_i$. 
Therefore, if one knows 
the asymptotic for  each $N(V_i, {\cal L}, B)$
then one immediately 
obtains the asymptotic for $N(V, {\cal L}, B)$.

\begin{dfn}
{\rm Let $V$ be an ${\cal L}$-primitive algebraic variety over 
a number field $F$, 
$\rho\,:\, X \ra \overline{V}^{\cal L}$ a desingularization  over $F$ 
of the  closure of $V$ with the exceptional locus consisting of 
smooth irreducible divisors 
$D_1,\ldots, D_l$. We consider ${\rm Pic}(X)$ and  ${\rm Pic}(V,{\cal L})$ as  
${\rm Gal}(\overline{F}/F)$-modules 
and we denote by $\beta_{\cal L}(V)$ 
the rank of ${\rm Gal}(\overline{F}/F)$-invariants 
in  ${\rm Pic}(V,{\cal L})$ and by 
$\delta_{\cal L}(V)$ the cardinality of the 
cohomology group 
\[ \H^1({\rm Gal}(\overline{F}/F),   
{\rm Pic}(V, {\cal L})). \]}
\label{kohom}
\end{dfn}

\begin{rem}
{\rm Using the long exact Galois-cohomology sequence associated with 
(\ref{sh-3}), one immediately obtains that $\beta_{\cal L}(V)$ 
and  $\delta_{\cal L}(V)$ 
do not  depend on the choice of the resolution $\rho$. 
} 
\end{rem}          

\noindent
{\bf Step 4 (expected asymptotic formula): } 
Let  $V$ be a strongly ${\cal L}$-saturated (and  ${\cal L}$-primitive) 
smooth quasi-projective variety. Assume that  $a_{\cal L}(V) >0$. 
Then we expect that the following asymptotic formula holds: 
$$
N(V, {\cal L}, B) = c_{\cal L}(V)B^{\alpha_{\cal L}(V)} 
(\log B)^{ \beta_{\cal L}(V)-1} \left( 1 + o(1) \right),   
$$
where  
$$
c_{\cal L}(V) := \frac{\gamma_{\cal L}(V)}{
 \alpha_{\cal L}(V)(\beta_{\cal L}(V)-1)!}
\delta_{\cal L}(V)\tau_{\cal L}(V),
$$
$\gamma_{\cal L}(V)$ is an invariant of the triple
$({\rm Pic}(V,{\cal L}), {\rm Pic}(V,{\cal L})_{\R},
\Lambda_{\rm eff}(V,{\cal L}))$ (\ref{gamma-dfn}), 
$\delta_{\cal L}(V)$ is a cohomological invariant of the
${\rm Gal}(\overline{F}/F)$-module 
${\rm Pic}(V, {\cal L})$ (\ref{kohom}) and 
$\tau_{\cal L}(V)$  is an adelic 
invariant of a family of $v$-adic 
metrics $\{ \|\cdot \|_v \}$ on $L$ (\ref{tau-dfn}). 
\medskip

In sections 4 and 5  we discuss some examples which show
how the constants 
$$
\a_{\cal L}(V), \b_{\cal L}(V),\delta_{\cal L}(V),
\gamma_{\cal L}(V),\tau_{\cal L}(V)
$$
appear in asymptotic formulas for the number of 
rational points of bounded ${\cal L}$-height
on algebraic varieties. Naturally, we expect 
that the exhibited behavior is typical. 
However, we also feel that one should collect 
more examples which could help
to clarify the general situation.

\subsection{${\cal L}$-primitive 
fibrations and $\tau_{\cal L}(V)$}

We proceed to discuss our 
observations concerning the arithmetic 
conjecture \ref{aaf} 
at its relation to the geometric 
conjecture \ref{conj-fb}.

Let $V$ be a weakly ${\cal L}$-saturated 
smooth quasi-projective 
variety with $a_{\cal L}(V) >0$ which is not strongly saturated
and which doesn't contain Zariski open dense strongly saturated
subvarieties. We distinguish the following 
two cases: 

\noindent
{\bf Case 1.  $V$ is not ${\cal L}$-primitive. }
In this case we expect that some Zariski open dense subset 
$U \subset V$ admits an 
${\cal L}$-primitive fibration which is defined by 
a projective regular 
morphism $\pi\, : \, X \ra Y$ over $F$ 
to a low-dimensional 
normal irreducible projective variety $Y$ 
satisfying the conditions 
(i)-(iii) in \ref{prim-fb}, for an appropriate 
smooth projective compactification 
$X$ of $U$ (see \ref{conj-fb}). 
It seems natural to expect that all fibers satisfy the vanishing 
assumption \ref{vanish}. 
Thus we see that for any $y \in Y(F)$ such that 
$V_y = \pi^{-1}(y) \cap V$ is 
${\cal L}$-primitive we can define
the adelic number $\tau_{\cal L}(V_y)$. 
Furthermore, we expect 
that every ${\cal L}$-target $W$ is contained 
in an appropriate 
${\cal L}$-primitive subvariety $V_y$ which is a fiber 
of the   ${\cal L}$-primitive fibration 
$\pi \; :\; V \ra U$ on $V$. 
In particular, Step 4 of our main strategy implies that if  every 
${\cal L}$-target $W$ {\em coincides} with  a
suitable ${\cal L}$-primitive fiber $V_y$ 
then  one should expect the asymptotic 
$$
N(V, {\cal L}, B) = c_{\cal L}(V)B^{\alpha_{\cal L}(V)} 
(\log B)^{ \beta_{\cal L}(V)-1} \left( 1 + o(1) \right),   
$$
where  the numbers $\alpha_{\cal L}(V)$ 
(resp.  $\beta_{\cal L}(V)$) 
coincide with the numbers $\alpha_{\cal L}(V_y)$ 
(resp.  $\beta_{\cal L}(V_y)$) for the corresponding ${\cal L}$-targets
$V_y$ and the constant $c_{\cal L}(V)$ is equal to the sum 
$$
\sum_{y}  c_{\cal L}(V_y) = 
\sum_{y}  \frac{\gamma_{\cal L}(V_y)}{
 \alpha_{\cal L}(V_y)(\beta_{\cal L}(V_y)-1)!}
\delta_{\cal L}(V_y)\tau_{\cal L}(V_y),
$$
where $y$ runs over all points in $Y(F)$ such that  $V_y$ is an 
${\cal L}$-target of $V$. 
It is natural to try to understand 
the dependence of $\tau_{\cal L}(V_y)$ on the choice of a point 
$y \in Y(F)$. We expect that the number 
$\tau_{\cal L}(V_y)$ can be 
interpreted as a ``height'' of $y$. 
More precisely, the examples we considered
suggest the  following:

\begin{conj}
There exist a family of $v$-adic metrics on $K_Y$ and two positive 
constants $c_2 > c_1 > 0$ such that 
\[      c_1 H_{\cal F}(y) \leq \tau_{\cal L}(V_y) \leq c_2 
H_{\cal F}(y) \;\; \forall y \in Y(F) \cap U, \]
where $U \subset Y$ is some dense Zariski open subset and 
${\cal F}$ is a metrized ${\bf Q}$-invertible sheaf 
associated with the ${\bf Q}$-Cartier divisor 
$L_1^{-1} \otimes K_Y$  
(recall that $L_1$ is the tautological  ample  
${\bf Q}$-Cartier divisor 
on $Y$ defined by the graded 
ring ${\rm R}(V, {\cal L})$ $($see 
\ref{proj-fb}$))$. 
\end{conj}

\noindent
{\bf Case 2.  $V$ is ${\cal L}$-primitive 
(but not strongly saturated!).} 
We don't know examples of a precise 
asymptotic formula in this situation.

\begin{exam}
{\rm Let $V$ be  a Fano diagonal cubic  bundle over 
${\bf P}^3$ with the homogeneous coordinates $(X_0:X_1:X_2:X_3)$ 
defined as a hypersurface in ${\bf P}^3 \times {\bf P}^3$ 
by the equation 
\[ X_0Y_0^3 + X_1Y_1^3 + X_2Y_2^3+ X_3Y_3^3= 0 \] 
in ${\bf P}^3 \times {\bf P}^3$ (see \cite{BaTschi4}). 
We expect that $V$   
is not strongly saturated 
with respect to a metrized  
anticanonical sheaf  $L:= {\cal O}(3,1)$ 
and that the corresponding 
${\cal L}$-targets are the splitting diagonal 
cubics in fibers of 
the natural projection $\pi\,: \, V \to {\bf P}^3$
(this leads to the failure of the expected
asymptotic formula in Step 4 for this example). 
}
\end{exam}

The next example  
was suggested to us by Colliot-Th\'el\`ene: 

\begin{exam}
{\rm Let $V$ be an analogous diagonal quadric  bundle over 
${\bf P}^3$ 
defined as a hypersurface in ${\bf P}^3 \times {\bf P}^3$ 
by the equation 
\[ X_0Y_0^2 + X_1Y_1^2 + X_2Y_2^2 + X_3Y_3^2 = 0. \] 
For infinitely many fibers $V_x = \pi_1^{-1}(x)$ ($x \in {\bf P}^3(F)$) 
we have ${\rm rk}\,{\rm Pic}(V_x) = 2$. At the same time, we have 
also ${\rm rk}\,{\rm Pic}(V) = 2$. We consider the height 
function associated to some metrization of the line bundle
$L:= {\cal O}(3,2)$.
On the one hand, we think that the asymptotic on the whole variety
is $c(V)B\log B(1+o(1))$ for $B\ra \infty$ with some $c(V)>0$. 
On the other hand, if $X_0 X_1 X_2 X_3$ is a square in $F$
we get already about $B\log B$ solutions. 
Another important observation
is the {\em expected} convergency of the series
\[ \sum_{x \in {\bf P}^3(F) : V_x \cong {\bf P}^1 \times {\bf P}^1} 
c(V_x). \]
The latter would be a consequence of the following two facts. 
First, the condition  $V_x \cong {\bf P}^1 \times {\bf P}^1$
($x = (X_0:X_1:X_2:X_3)$) is equivalent to the conditions that $V_x$ 
contains an $F$-rational point and that the product 
$X_0 X_1 X_2 X_3$ is a square in $F$.  The number of $F$-rational points 
$x'= (X_0:X_1:X_2:X_3:Z)$  with $H_{{\cal O}(1)}(x') \leq B$ 
lying on the hypersurface with the 
equation $X_0X_1X_2X_3 = Z^2$ in the weighted projective space 
${\bf P}^4(1,1,1,1,2)$ can be estimated from above by $B^2 
(\log B)^3( 1 + o(1))$. 

\noindent
Secondly, we expect  
$$
c(V_x) = H^{-1}_{{\cal O}(4)}(x) 
( 1 + o(1) ), 
$$ 
{\em uniformly} over the base ${\bf P}^3(F)$. 
This would imply the claimed convergency. 
}
\end{exam}


\section{Height zeta-functions}

\subsection{Tauberian theorem}

One of the main techniques in the proofs of asymptotic formulas 
for the counting function
$$
N(V,{\cal L},B):= \#\{x\in V(F)\,\,:\,\, H_{\cal L}(x)\le B\,\}
$$
has been the use of {\em height zeta functions}. Let 
${\cal L}$ be an ample metrized invertible sheaf 
on a smooth quasi-projective
algebraic variety $X$. We define the height zeta function
by the series
$$
Z(X,{\cal L},s):=\sum_{x\in X(F)}H_{\cal L}(x)^{-s}
$$
which converges absolutely for ${\rm Re}(s)\gg 0$. 
After establishing the analytic properties  
of $Z(X,{\cal L},s)$ one uses the
following version of a Tauberian theorem:

\begin{theo}(\cite{delange})  Suppose that 
there exist an $\e >0$ and a real number
$\Theta({\cal L})>0$ such that
$$
Z(X,{\cal L},s)=\frac{\Theta({\cal L})}{(s-a)^{b}}
+ \frac{f(s)}{(s-a)^{b-1}}
$$
for some $a>0$, $b\in \N$ and some function  $f(s)$ 
which is holomorphic for ${\rm Re}(s)>a-\e$.
Then we have the following asymptotic formula
$$
N(X,{\cal L}, B)= \frac{\Theta({\cal L})}{a\cdot (b-1)!}B^a(\log B)^{b-1}(1+o(1))
$$
for $B\ra \infty$. 
\end{theo}

\subsection{Products}

Let  $X_1$ and $X_2$ be two smooth quasi-projective 
varieties with ample metrized invertible sheaves 
${\cal L}_1$ and ${\cal L}_2$ (resp. on $X_1$ and $X_2$). 
Denote by $X=X_1\times X_2$ the product and by ${\cal L}$ 
the product of ${\cal L}_1$ and ${\cal L}_2$ (with the obvious 
metrization). Clearly,
\[  Z(X,{\cal L}, s)=  Z(X_1,{\cal L}_1, s)\cdot Z(X_2,{\cal L}_2, s). \]
Assume that for $i=1,2$ we have 
$$
Z(X_i,{\cal L}_i,s)= \frac{\Theta({\cal L}_i)}{
(s-\alpha_{{\cal L}_i}(X_i))^{\b_{{\cal L}_i}(X_i)}}
+ \frac{f_i(s)}{(s-\alpha_{{\cal L}_i}(X_i))^{\b_{{\cal L}_i}(X_i)}}
$$
with 
some functions $f_i(s)$ which are holomorphic in the domains
${\rm Re}(s_i)>\alpha_{{\cal L}_i}(X_i)-\e$
for some $\e>0$.  There are two possibilities: 

\smallskip
Case 1: $\alpha_{{\cal L}_1}(X) = \alpha_{{\cal L}_2}(X)$.
In this situation the constant $\Theta({\cal L})$ at the pole of 
highest order $\b_{{\cal L}_1}(X_1)+\b_{{\cal L}_2}(X_2)$ is 
given by $\Theta({\cal L})=\Theta({\cal L}_1)\Theta({\cal L}_2)$.

\medskip
Case 2: 
$\alpha_{{\cal L}_1}(X) <  \alpha_{{\cal L}_2}(X)$. In this situation 
the constant is a sum  
$$ 
\Theta({\cal L}) = \sum_{x \in X_1(F)}
H_{{\cal L}_1}^{-\alpha_{{\cal L}_2}}(x)\Theta({\cal L}_2). 
$$
Consider the projection $X\ra X_1$ and denote by $V_x$ the fiber over
$x\in X_1(F)$. We notice
$$
\tau_{\cal L}(V_x) = H_{{\cal L}_1}^{-\alpha_{{\cal L}_2}}(x)
\tau_{{\cal L}_2}(X_2).
$$
We denote by $M$ the ${\bf Q}$-Cartier divisor 
$\pi_1^*{\cal L}_1^{\alpha_{{\cal L}_2}} \otimes K_{V_1}$. 
We obtain that 
$\pi_1^*{\cal L}_1^{\alpha_{{\cal L}_2}} = M \otimes  K_{V_1}^{-1}$. 
So we have 
$\tau_{\cal L}(V_x) \sim H_{M  \otimes  K_{V_1}^{-1}}^{-1}$.
We observe that Tamagawa numbers of fibers  
depend on the height of the points on the base.

\subsection{Symmetric product of a curve}

Let $C$ be a smooth irreducible curve of genus $g\ge 2$ over $F$. 
We denote by $X=C^{(m)}$ the $m$-th symmetric product of $C$
and by  $Y:={\rm Jac}(C)$ the Jacobian of $C$. 
We fix an $m>2g-2$. We have a fibration 
$$
\pi\, :\; C^{(m)} \,\ra \, Y,
$$
with $\P^{m-g}$ as fibers. We denote by $V_y$ a fiber over $y\in Y(F)$. 

Let $\tilde{C}\ra {\rm Spec}({\cal O}_F)$  be a smooth
model of $C$ over the integers and ${\cal L}$  an ample hermitian 
line bundle on $\tilde{C}$. It defines a height function
$$
H_{\cal L}\,:\,C(F)\ra \R_{>0}
$$ 
which extends to 
$X(F)$.  
Observe that $\a_{\cal L}(X)=(m+1-g)/d$, where  
$d:={\rm deg}_{\Q}({\cal L})$. 
Consider the height 
zeta function
$$
Z(X,{\cal L},s):=\sum_{x\in X(F)}H_{\cal L}(x)^{-s}.
$$
This function was  introduced by Arakelov  in \cite{arakelov}.

\begin{theo}\cite{faltings}
Let ${\cal L}$ be an ample hermitian line bundle on $\tilde{C}$. 
There exist an $\e ({\cal L})>0$ and a real number
$\Theta({\cal L})\neq 0$ such that the height zeta function 
has the following representation
$$
Z(X,{\cal L},s)=\frac{\Theta({\cal L})}{
(s-\a_{\cal L}(X))} +f(s)
$$
with some function $f(s)$ which is holomorphic 
for ${\rm Re}(s)>\a_{\cal L}(X)-\e({\cal L})$.

\end{theo}

This is Theorem 8 in (\cite{faltings}, p. 422). 
Arakelov gives an explicit expression for the constant
$\Theta({\cal L})$ (\cite{arakelov}). We are very grateful to
J.-B. Bost for pointing out to us that Arakelov's formula is 
not correct and for allowing us to use his notes on the
Arakelov zeta function. 

\begin{theo} With the notations above we have
$$
\Theta({\cal L})=\sum_{y\in Y(F)}\tau_{\cal L}(V_y).
$$
\end{theo}

\noindent
{\em Proof.} We outline the proof for 
$F=\Q$. For ${\rm Re}(s)\gg 0$ one 
can rearrange the order of summation
and one  obtains
$$
Z(X,{\cal L},s):=\sum_{y\in Y(F)}
\sum_{x\in V_y(F)}H_{\cal L}(x)^{-s}.
$$
It is proved in (\cite{faltings}, p. 420-422)
that the sums
$$
Z(V_y,{\cal L},s):=\sum_{x\in V_y(F)}H_{\cal L}(x)^{-s}
$$
have  simple poles at $s=\a_{\cal L}(X)$
with  non-zero residues 
and that one can ``sum'' these expressions to 
obtain a function with a simple pole at $s=\a_{\cal L}(X)$ and
meromorphic continuation to ${\rm Re}(s)>\a_{\cal L}(X)-\e({\cal L})$
for some $\e({\cal L})>0$.
Moreover,  the residue at this pole is
obtained as a sum over $y\in Y(F)$ of the residues of 
$Z(V_y,{\cal L},s)$. 

Choosing an element $z(y)$ in the class of $y\in {\rm Jac}(C)(F)$,
and denoting by $E(y):=\Gamma (C,{\cal O}(z(y)))$ 
one can identify the fiber as
$V_y= \P(E(y))$,
where $f\in E(y)= \Gamma (C,{\cal O}(z(y)))$ is mapped to 
$z(y)+{\rm div}(f)$. 
The height is given by the formula
$$
H_{\cal L}(z(y)+{\rm div}(f)):= H_{\cal L}(z(y))
\exp(\int_{C(\C)}\log |f|c_1({\cal L})).
$$
This defines a metrization of the anticanonical 
line bundle on $V_y= \P(E(y))$. Assuming the Tamagawa number 
conjecture  for $\P^n$ (for suitable metrizations
of the line bundle ${\cal O}(1)$ on $\P^n$) we obtain   
$$
\lim_{s\ra \a_{\cal L}(X)}(s-\a_{\cal L}(X))Z(V_y,{\cal L},s)=
\tau_{\cal L}(V_y).
$$
\hfill $\Box$

\medskip
\noindent
One can write down an explicit formula for
$\tau_{\cal L}(V_y)$.
For $f\in E(y)_{\C}\backslash \{0\}$ 
we define 
$$
\Phi(f):=\exp(\frac{1}{d}\int_{C(\C)}\log |f|c_1({\cal L}))
$$
and we put $\Phi(0)=0$. 
It follows that
$$
\tau_{\cal L}(V_y)=
\frac{1}{2}\a_{\cal L}(X)H_{\cal L}(y)^{-\a_{\cal L}(X)}
\cdot 
\frac{{\rm vol}(\{ f\in E(y)_{\R} \,|\, \Phi(f)\le 1\})}{
{\rm vol}(E(y)_{\R}/ E(y))},
$$
where the volumes are calculated with respect to some 
Lebesgue measure on $E(z(y))_{\R}$. 
Arakelov relates this last expression to the Neron-Tate height of
$y\in {\rm Jac}(C)$. 
A detailed calculation due to Bost indicates that Arakelov's formula
is correct only up-to $O(1)$.

\subsection{Homogeneous spaces $G/P$}

Let $G$ be a split semisimple linear algebraic group defined over a 
number field $F$. It contains a Borel subgroup 
$P_0$ defined over $F$ and a maximal torus which is 
split over $F$. Let $P$ be a standard parabolic. 
Denote by $Y_P=P\backslash G$ (resp. $X=P_0\backslash G $) 
the corresponding flag variety. 

A choice of a 
maximal compact subgroup ${\bf K}$ such that 
$G({\bf A}_F)=P_0({\bf A}_F){\bf K}$ defines a metrization on 
every line bundle $L$ on the flag varieties $Y_P$ (\cite{FMT}, p. 426).
We will denote by 
$$
H_{\cal L}\,:\, P(F)\backslash G(F)\ra \R_{>0}
$$ 
the associated height. We consider the height zeta function
$$
Z(X,{\cal L},s):= 
\sum_{x \in X(F)} H_{\cal L}(x)^{-s}.
$$

\begin{theo} Let ${\cal L}$ be an ample metrized line bundle on $X$.
There
exist an $\e ({\cal L})>0$ and a real number
$\Theta({\cal L})\neq 0$ such that the height zeta function 
has the following representation
$$
Z(X,{\cal L},s)=\frac{\Theta({\cal L})}{(s-
\a_{\cal L}(X))^{\b_{\cal L}(X)}} +
\frac{f(s)}{(s-\a_{\cal L}(X))^{\b_{\cal L}(X)-1}}
$$
with some function $f(s)$ which is holomorphic 
for ${\rm Re}(s)>\a_{\cal L}(X)-\e({\cal L})$.
\end{theo}

This theorem follows from 
the identification of the height zeta function
with an Eisenstein series and from the work of Langlands. 
The formula (2.10) in (\cite{FMT}, p. 431) provides
an expression for $\Theta({\cal L})$ which we will now
analyze.  

There is a canonical way to identify the 
faces of the closed cone of effective
divisors $\L_{\rm eff}(X)\subset 
{\rm Pic}(X)_{\R}$  
with $\L_{\rm eff}(Y_P)$ as $P$ runs through
the set of standard parabolics. 
A line bundle $L$ such that its class is
contained in the interior of the cone 
$\in \L_{\rm eff}(X)$ defines a line
bundle 
$$
[L_Y]:=\a_{\cal L}(X)[L]+ [K_{X}]
$$
which is contained in the interior of the
face $\L_{\rm eff}(Y)\subset 
\L_{\rm eff}(X)$ for some $Y=Y_P$.  
We have a fibration
$\pi_{\cal L}\,: X\ra Y$ 
with fibers isomorphic to the flag variety $V:=P_0\backslash P$. 
A fiber over $y\in Y(F)$ will be denoted by $V_y$. 
Denote by ${\cal K}_Y$ the canonical line bundle on 
$Y$ with the metrization defined above.

\begin{theo}
We have 
$$
\Theta({\cal L})= \sum_{y\in Y(F)}\gamma_{\cal L}(X)\tau_{\cal L}(V_y).
$$
\end{theo}

\noindent
{\em Proof.}
In the domains of absolute
and uniform convergence we can rearrange the order of summation and
we obtain
$$
Z(X,{\cal L},s)=\sum_{y\in Y(F)}
\sum_{x\in P_0(F)\backslash P(F)} H_{\cal L}(yx)^{-s}.
$$
One can check that the sums
$$
Z(V_y,{\cal L},s):=\sum_{x\in V_y(F)}H_{\cal L}(x)^{-s}
=\sum_{x\in P_0(F)\backslash P(F)} H_{\cal L}(yx)^{-s}
$$
have poles at $s=\a_{\cal L}(X)$ of order $\b_{\cal L}(X)$
with  non-zero residues, and that they admit meromorphic 
continuation to  $\a_{\cal L}(X)-\e({\cal L})$ for some
$\e({\cal L})>0$. Moreover,
the constant $\Theta({\cal L})$ is
obtained as a sum over $y\in Y(F)$ of the residues of 
$Z(V_y,{\cal L},s)$. From the Tamagawa number 
conjecture  for $P_{0}\backslash P$ (with varying metrizations
of the anticanonical line bundle) we obtain   
$$
\lim_{s\ra \a_{\cal L}(X)}
(s-\a_{\cal L}(X))^{\b_{\cal L}(X)}Z(V_y,{\cal L},s)=\gamma_{\cal L}(X)
\tau_{\cal L}(V_y).
$$
The sum 
$$
\sum_{y\in Y(F)}\gamma_{\cal L}(X)
\tau_{\cal L}(V_y)
$$
converges for $[L_Y]$ contained in the interior of 
$\L_{\rm eff}(Y)$. 
\hfill $\Box$
\medskip

\noindent
Let us recall the explicit formula for $\tau_{\cal L}(V_y)$ 
(see (2.9) in \cite{FMT}, p. 431):
$$
\lim_{s\ra \a_{\cal L}(X)}(s-\a_{\cal L}(X))^{
\b_{\cal L}(X)} \sum_{x\in P_0(F)\backslash P(F)}
H_{\cal L}(xy)^{-s}
=\g_1 c_P H_{{\cal L}_Y^{-1}\otimes {\cal K}_{Y}}(y)
$$
where $\g_1\in \Q$
is an explicit constant
and the constant $c_P$ is defined in (\cite{FMT}, p. 430).
It follows that
$$
\Theta({\cal L})=
\g_1 c_P\sum_{y\in Y(F)}
H_{{\cal L}_Y^{-1}\otimes {\cal K}_Y}(y).
$$
Next we observe that there
is an explicit constant $\g_2\in \Q$ such that we have 
$$
\g_2\tau_{\cal L}(V_{y})=
c_P\cdot H_{{\cal L}^{-1}_Y\otimes {\cal K}_Y}(y)
$$
for all $y\in Y(F)$. 
To see this, we first identify 
$c_P= \g_2\tau_{\cal L}(V)$, this is done by a computation 
of local factors of intertwining operators (\cite{peyre}, p.160-161).
The next step involves the comparison of Tamagawa measures on 
$V_{y}$ for varying $y\in Y(F)$. 
Finally, we have $\g_{\cal L}(X)=\g_1\g_2$.

\subsection{Toric varieties}

There are many equivalent ways to describe a toric variety $X$
over a number field $F$ together with
some projective embedding (see, for example,  \cite{cox,fulton,danilov}). 
For us, it will be useful to view
$X=X_{\S}$ as 
a collection of the following data (\cite{BaTschi1}): 

1. A splitting field $E$ of the  algebraic torus $T$ and the Galois group
$G= {\rm Gal} (E/F)$.

2. The lattice of $E$-rational characters of $T$, which we denote by $M$
and its dual lattice $N$.

3. A $G$-invariant complete
fan $\S$ in $N_{\R}=N\otimes \R$.

There is an isomorphism between the group of 
$G$-invariant integral piecewise linear functions $\varphi\in PL(\S)^G$ 
and classes of $T$-linearized line bundles on $X_{\S}$. 
For $\varphi\in PL(\S)^G$ 
we denote the corresponding line bundle by $L(\varphi)$.

We define metrizations of line bundles as follows.
Let $G_v\subset G$ be the decomposition group at $v$. 
We put $N_v=N^{G_v}$ 
for the lattice of $G_v$-invariants of $N$ 
for non-archimedean valuations $v$ 
and  $N_v=N^{G_v}_{\R}$ for archimedean $v$. 
We have the logarithmic map
$$
T(F_v)/T({\cal O}_v)\ra N_v
$$
which is an embedding of finite index for 
all non-archimedean $v$, an isomorphism of lattices for
almost all non-archimedean valuations 
and an isomorphism of real vector spaces for archimedean valuations.  
We denote by $\overline{t}_v$ the image of $t_v\in T(F_v)$ under
this map.

\begin{dfn}(\cite{BaTschi1}, p. 607)
{\rm 
For every $\varphi\in PL(\S)^G$  and 
$t_v\in T(F_v)$ we define the local height 
function 
$$
H_{\S,v}(t_v,\varphi):=e^{\varphi(\overline{t}_v)\log q_v}
$$
where $q_v$ is the cardinality of the residue field of $F_v$ for
non-archimedean valuations and $\log q_v=1 $  for archimedean  
valuations. For $t\in T({\bf A}_F)$ we define
the global height function as
$$
H_{\S}(t,\varphi):=\prod_{v\in {\rm Val}(F)}H_{\S,v}(t_v,\varphi).
$$
}
\end{dfn}

\noindent
We proved in (\cite{BaTschi1}, p. 608) that this pairing
can be extended to a pairing
$$
H_{\S}\,:\, T({\bf A}_F)\times PL(\S)^G_{\C}\ra \C
$$
and  that  it
defines a simultaneous metrization of $T$-linearized line bundles
on $X$. We will denote such metrized line bundles 
by ${\cal L}={\cal L}(\varphi)$.
We consider the height zeta function
$$
Z(T,{\cal L},s)=\sum_{t\in T(F)}H_{\cal L}(t)^{-s}.
$$

\begin{theo}(\cite{BaTschi3})
Let ${\cal L}$ be an invertible sheaf on $X$
(with the metrization introduced above) such that its class  
$[L]$ is contained in the interior of $ \L_{\rm eff}(X)$.
There exist an $\e({\cal L})>0$ and a $\Theta({\cal L})>0$ 
such that the height zeta function 
has the following representation
$$
Z({\cal L},T,S)=\frac{\Theta({\cal L})}{
(s-\a_{\cal L}(X))^{\b_{\cal L}(X)}}+
\frac{f(s)}{(s-\a_{\cal L}(X))^{\b_{\cal L}(X)-1}}
$$
where $f(s)$ is a function which is holomorphic for
${\rm Re}(s)>\a_{\cal L}(X)-\e({\cal L})$. 
\end{theo}

\begin{rem}{\rm
The computation of the constants
$\a_{\cal L}(X)$ and $\b_{\cal L}(X)$ 
in specific examples is
a problem in linear programming. 
For the anticanonical line
bundle on a smooth toric variety $X$ we have 
$\a_{{\cal K}^{-1}}(X) =1$ and 
$\b_{{\cal K}^{-1}}(X) =
\dim PL(\S)^G_{\R}-\dim M^G_{\R}$.
}
\end{rem}

Our goal is to identify the constant $\Theta({\cal L})$.
Let us recall some properties of toric
varieties and introduce more notations (see \cite{BaTschi1}). 
The cone of effective divisors $\L_{\rm eff}(X)$ is generated by 
the classes of irreducible components of $X\backslash T$ which we
denote by $[D_1],...,[D_r]$. These divisors correspond to
Galois orbits $\S_1(1),...,\S_r(1)$ 
on the set of $1$-dimensional cones in $\S$.
The line bundle ${\cal L}$ defines a face $\L({\cal L})$ of $\L_{\rm eff}(X)$.
We denote by $J=J({\cal L})$ the 
maximal set of indices $J\subset [1,...,r]$
such that we have
$$ 
\a_{\cal L}(X)[L]+[K_X]=\sum_{j\in J}r_j[D_j]
$$
with  $[D_j]\in \L({\cal L})$ and some $r_j\in \Q_{>0}$. 
We denote by $I=I({\cal L})$ the set of indices
$i\not\in J({\cal L})$. 
We denote by $M_J$ the lattice given by
$$
M_J:=\{ m\in M\,|\, <e,m>=0\hskip 0,3cm 
\R_{\ge 0}e\in \cup_{i\in I } \S_i(1) \},
$$
We denote by $M_I:=M/M_J$ and by $N_*$ the corresponding dual
lattices. We have an exact sequence  of
algebraic tori
$$
1\ra T_I\ra T\ra T_J\ra 1 
$$
which induces a map
$\pi_{\cal L}\,:\, T(F)\ra T_J(F)$
with finite cokernel 
and an exact sequence of lattices
$$
0\ra N_I\ra N\ra N_J\ra 0.
$$ 
The restriction of the fan $\S\subset N_{\R}$ to 
$N_{I,\R}$ will be denoted by $\S_I$. It is again 
a $G$-invariant fan and it
will define an equivariant compactification  $X_I$ of $T_I$.
The class of the piecewise linear function $\varphi_I\in PL(\S_I)^G$ with
$\varphi_I(e)=1$ for $e\in \cup_{i\in I}\S_i$ corresponds to the class of
the anticanonical line bundle $[-K_I]\in {\rm Pic}(X_I)$. 

The line bundle ${\cal L}$ defines a fibration  of 
varieties $\pi_{\cal L}\,:\, X_{\S}\ra Y$ with fibers isomorphic to $X_I$,
which, when restricted to $T$,
gives rise to the exact sequence of tori above. We denote the
fiber over $y\in T_J(F)$  by $X_{I,y}$.

\begin{theo} We have
$$
\Theta({\cal L})=\g_{\cal L}(X)
\d_{\cal L}(X)\sum_{y\in \pi_{\cal L}(T(F))}\tau_{\cal L}(X_{I,y}).
$$
\end{theo}

{\em Proof.} Let ${\cal L}={\cal L}(\varphi)$ be
an invertible sheaf on $T$ with the metrization introduced above.
In the domain of absolute and uniform convergence
we can rearrange the order of summation 
and we obtain
$$
Z(T,{\cal L},s)=\sum_{y\in \pi_{\cal L}(T(F))}
\sum_{x\in T_{I,y}(F)} H_{\S}(yx,\varphi)^{-s}.
$$
From the proof of our main theorem in \cite{BaTschi3} 
it follows that the sums
$$
Z(T_{I,y},{\cal L},s)=\sum_{x\in T_{I}(F)} H_{\S}(yx,\varphi)^{-s}
$$
have a pole at $\a_{\cal L}(X)$ of order 
$\b_{\cal L}(X)$. 
Moreover, the constant $\Theta({\cal L})$ is obtained
as a  sum over $y\in \pi_{\cal L}(T(F))$ of residues of 
$Z(T_{I,y},{\cal L},s)$. 
Now we want to use the Tamagawa number
conjecture for the anticanonical line bundle 
(with varying metrizations) on the
toric variety $X_I$ to conclude the proof. 

In \cite{BaTschi2} we proved this conjecture
for a specific metrization and under 
the assumption that the fan $\S$ is {\em regular}.
We want to demonstrate that our proof goes through in the general case
needed above. 

\smallskip
Our main idea was to use the Poisson summation formula on the
adelic group $T({\bf A}_F)$ and to obtain an integral 
representations for the height zeta function.
We denote by 
${\cal A}_I=(T_I({\bf A}_F)/{\bf K}T_I(F))^*$ the group of
unitary characters of $T_I({\bf A}_F)$ which are
trivial on $T_I(F)$ and on the maximal compact subgroup 
${\bf K}\subset T_I({\bf A}_F)$. Using the
the adelic definition of the height
function we obtain 
$$
Z(T_{I,y},{\cal L},s)=\sum_{t\in T_I(F)}H_{\cal L}(yt)^{-s}
=\int_{{\cal A}_I}d\chi \int_{T_I({\bf A}_F)} H_{\cal L}(yt)^{-s}\chi(t)d\mu,
$$
where $d\mu$ is a Haar measure
on $T({\bf A}_F)$ and $d\chi$ is the orthogonal 
Haar measure on ${\cal A}_I$. 
To apply our technical theorem in \cite{BaTschi2} about the
analytic continuation and the residues of such integrals
we need to know that
$$
\int_{T_I({\bf A}_F)} H_{\cal L}(yt)^{-s}\chi(t)d\mu =\prod_{i\in I}
L_i(\chi_i,s) \cdot \zeta_{\S_I}(\chi,s)\cdot \zeta_{\infty}(\chi,s)
$$
where 
$$
\zeta_{\S_I}(\chi,s)=\prod_{v\in {\rm Val}(F)}\zeta_{\S_I,v}(\chi,s)
$$ 
is an absolutely convergent 
Euler product for ${\rm Re}(s)>\a_{\cal L(X)}-\e({\cal L})$,
$L_i(\chi_i,s)$ are Artin $L$-functions 
(with some induced characters $\chi_i$)
and $\zeta_{\infty}(\chi,s)$ satisfies certain growth conditions. 

First we observe that it is unnecessary to assume that 
the fan $\S_I$ is regular. Using our definition of the height function
we see that the calculation of the Fourier transform (see \cite{BaTschi1})
reduces to summations of the function 
$q_v^{{\varphi}(\overline{t_v}) +im_v}$
($m_v\in M_{I,\R}$) over the lattice $N_{I,v}$
(resp. to integrations in cones
for archimedean valuation). A piecewise linear function
$\varphi$ induces a piecewise linear function 
on any subdivision of the fan. Clearly, the result of such summations
and integrations does not depend on any subdivisions.

Next we see that for a fixed $y\in T_J(F)$ we have
$H_{{\cal L},v}(yt)=H_{{\cal L},v}(t)$ for almost 
all $v$ and all $t\in T_I({\bf A}_F)$. 
Now we can refer to lemma 5.10 in \cite{BaTschi3} which proves
the required statement. 
The local integrals
for the remaining finitely many non-archimedean valuations
will be absorbed into $\zeta_{\S_I}(\chi,s)$. 
And finally, we need to check that the estimates 
of the Fourier transform of $H_{{\cal L},v}(yt)$ 
at archimedean valuations are still satisfied for any $y\in T_J(F)$.
This is straightforward.  

We  can now apply the main technical theorem of \cite{BaTschi2}
and obtain
$$
\lim_{s\ra \a_{\cal L}}(s-\a_{\cal L}(X))^{
\b_{\cal L}(X)}Z(T_{I,y},{\cal L},s)
= \gamma_{\cal L}(X) \delta_{\cal L}(X)\tau_{\cal L}(X_{I,y}).
$$ 
\hfill $\Box$

\begin{rem}
{\rm 
It is possible to compute $\tau_{\cal L}(X_{I,y})$ and to observe
that it is related to the  height of the point $y\in T_J(F)$. 
}
\end{rem}

\begin{rem}
{\rm
Similar statements hold for equivariant compactifications 
of homogeneous spaces $G/U$ where $G$ is a split 
reductive group and $U$ is its maximal 
unipotent subgroup \cite{StrTschi}. 
We hope that these results can
be extended to equivariant
compactifications of other homogeneous spaces, in particular,
to equivariant compactifications of reductive and non-reductive groups.
}
\end{rem}

\section{Singular Fano varieties}

\subsection{Weighted projective spaces}

Let $W:= {\bf P}(w) = {\bf P}(w_0, \ldots, w_n)$ be a weighted 
projective space of dimension $n$ with weights 
$w = (w_0, \ldots, w_n)$. We remark that $W$ is a rational variety 
over ${\bf Q}$ 
with ${\rm Pic}(W) \cong {\bf Z}$. Moreover, the anticanonical 
class $K_W^{-1}$ is an ample ${\bf Q}$-Cartier divisor. So $W$ 
is a (singular) Fano variety  of index
\[ r = \frac{ w_0 + \cdots + w_n} { l.c.m.\{ w_0, \ldots , w_n \}}. \]

One could try to generalize  the method of Schanuel \cite{schan} 
for counting ${\bf Q}$-rational points of bounded height on  
usual projective spaces to the case of weighted projective spaces. 
Let  $z_0, z_1, \ldots, z_n$ be homogeneous coordinates on 
$Y$. Then a first approximation to counting points of 
bounded height would be a counting of 
all $(n+1)$-tuples 
$(x_0, x_1, \ldots, x_n) \in {\bf Z}^n \setminus \{0 \}$ 
satisfying the conditions 
\[ |x_i| \leq 
B^{\frac{w_i}{w_0 + w_1 + \cdots + w_n}}\;\; i =0, \ldots, n. \]
Since the volume of the domain restricted by these inequalities is 
\[ B = \prod_{i =0}^n B^{\frac{w_i}{w_0 + w_1 + \cdots + w_n}},  \]
one  could expect that the asymptotic number 
of solutions of these inequalities 
agrees with ``expected'' 
linear growth for the anticanonical height. 
However, this ``intuition'' 
turns  out to be wrong, in general, because 
the singularities of $Y$ could be even worse 
than canonical. 
A typical class of singularities that appear on $Y$ are 
so called log-terminal 
singularities introduced by Kawamata \cite{Ka}.  
We give below a simple example of  a Del Pezzo surface 
with a log-terminal singularity 
and we show that for every dense Zariski open subsets $U \subset W$ the 
number $N(U,B)$ of $F$-rational points of  anticanonical height $\leq B$ 
in $U$ has more than linear growth: 
\[ N(U, B) = c(U)B^{2 - \frac{4}{m+2}}(1 + o(B)). \] 
Moreover, there are no 
dense Zariski open subsets 
$U' \subset X$ such that the adelic 
term in the constant $c(U)$ 
in the asymptotic formula for $N(U,B)$ 
would be independent of 
$U$ for $U \subset U'$. 

\begin{exam} 
{\rm {\sc (Del-Pezzo surface with a log-terminal singularity)} 
Let $W = {\bf P}(1,1,m)$ be a singular weighted projective plane 
with weights $(1,1,m)$, $m \geq 2$. Then the anticanonical class of 
$W$ is an ample ${\bf Q}$-Cartier divisor (i.e., $W$ is  Del Pezzo 
surface) and  $p = (0:0:1)$ is 
the unique singularity of $W$.  
Let $X \rightarrow W$ be the minimal resolution  
of the singularity at $p \in W$. Then $X$  is isomorphic 
to a ruled surface 
${\bf F}_m = {\bf P}({\cal O}_{{\bf P}^1} \oplus 
{\cal O}_{{\bf P}^1}(m))$ and the exceptional divisor $E = 
f^{-1}(p)$ is a smooth rational curve which is a section 
of the ${\bf P}^1$-bundle over ${\bf P}^1$ and  
$\langle E, E \rangle = -m$. Then we have 
\[ K_X = f^*K_W + \frac{2-m}{m} E. \]
Therefore,  $p$ is canonical $\Leftrightarrow$ 
$m = 2$ and $p$ is log-terminal  $\Leftrightarrow$ $m \geq 3$.

The group ${\rm Pic}(X)$ is isomorphic to ${\bf Z}^2$ where $\{ [E], 
[C] \}$ is a ${\bf Z}$-basis. Moreover,  $ [E], [C]$
are  generators for the cone $\Lambda_{\rm eff}(X)$ 
of effective divisors in ${\rm Pic}(X)_{\bf R}$. We have 
\[ K_X = -2E - (m+2 ) C, \]
\[ L := f^*(-K_X) = \frac{m + 2}{m}E +  (m + 2)C \]
and  
\[ {\alpha}_L(W) = \inf \{ t \in {\bf Q} \; : \; t[L] +[K_X] 
\in \Lambda_{\rm eff}(X) \} = 
\frac{2m}{m+2}.  \]
}
\end{exam}

Since $X$ is a smooth toric variety, we can apply our main result 
in \cite{BaTschi3} and obtain the following: 

\begin{theo}
Let $\pi: W \setminus p \ra {\bf P}^1$ be the natural projection, 
$C_x$ the fiber of $\pi$ over $x \in {\bf P}^1({\bf Q})$. Then   
for any dense Zariski open subset in $W \setminus 
p$, one has 
\[ N(U,B) = c(U)B^{2 - \frac{4}{m+2}}(1 + o(B)) \] 
Moreover,  
\[ c(U) =  \sum_{x \in {\bf P}^1({\bf Q})\cap \pi(U)}  c(C_x).   \]
\end{theo}

\subsection{Vaughan-Wooley cubic}

\begin{exam}
{\rm 
Let $Y \subset {\bf P}^5$ be a singular cubic defined by the 
equation $z_0z_1z_2 - z_3z_4z_5 = 0$, $X$ the intersection of 
$Y$ with the linear subspace in ${\bf P}^5$ with the equation: 
\[ z_0 + z_1 + z_2 - z_3 - z_4 - z_5 = 0. \]
Vaughan and Wooley proved \cite{VW}: 

\begin{theo}
Let $U \subset X$ be the complement in $X$ 
to the following $15$ divisors 
$D_{i_1 i_2 i_3}, D_{ij} 
\subset X$ $(\{ i_1, i_2, i_3 \} = \{ 0,1,2 \}, \; 
i \in \{ 0,1,2 \},\; j \in \{ 3,4,5\})$, where 
\[ D_{i_1 i_2 i_3} = \{ (z_0: \ldots : z_5) \in {\bf P}^5\; : \; 
z_{i_1} = z_3, 
\, z_{i_2} = z_4, \, z_{i_3} = z_5 \} \]
 
\[ D_{ij}   = \{ (z_0: \ldots : z_5) \in X\; : \; 
z_i = z_j = 0 \}.\]
If $N(U,B)$ is the number of ${\bf Q}$-rational points 
in $U$ of the anticanonical height $\leq B$, then 
there exist some  constants $c_1 > c_2 > 0$ such that 
\[   c_2 B^2 (\log B)^5 \leq  N(U,B) \leq c_1 B^2 (\log B)^5. \]
\end{theo}

We want to show that this result is compatible with predictions 
in \cite{BaMa}. First of all we note that $Y$ is a 
$4$-dimensional toric Fano  
variety: an equivariant compactification of a $4$-dimensional algebraic 
torus $T$ with respect to a   $4$-dimensional polyhedron $\Delta$ 
with $6$ lattice vertices 
\[ v_0 = (0,0,0,0), \,  v_1 = (1,0,0,0), \,   v_2 = (0,1,0,0), \] 
\[ v_3 = (0,0,1,0), \, v_4 = (0,0,0,1), \, v_5 = (1,1,-1,-1) \]
($\Delta$ is  the support of global 
sections of a very ample divisor $Y$ corresponding to the embedding  
$Y \hookrightarrow {\bf P}^4$). 
The polyhedron $\Delta$ has $9$ faces $\Theta_{ij}$ of codimension $1$: 
\[  \Theta_{ij} = {\rm Conv}(\{ v_0,v_1,v_2,v_3, v_4,v_5 \} \setminus 
\{v_i, v_j \} ). \]
Each face $\Theta_{ij}$ defines an torus invariant divisor $Y_{ij} \subset 
Y \setminus T$
such that $D_{ij} = Y_{ij} \cap X$.  
It is easy to check that all singularities of $Y$ are at worst terminal and 
that the hypersurface 
$X \subset Y$ intersects all strata $Y_{ij}$ transversally. 
From these facts we obtain that the only 
exceptional  divisors with the discrepancy $0$ 
that appear in a resolution of singularities of $X$ come from 
singularities in $X \cap T$. We write down the affine equation of  $X 
\cap T$ as 
\[ 1 + x + y - z -t - \frac{xy}{zt}, \]
where $T = {\rm Spec}\, {\bf Q}[x^{\pm 1}, y^{\pm 1},z^{\pm 1}, t^{\pm 1}]$. 
From this equation one immediately 
sees that  the only singularities 
in  $X \cap T$ are the $A_1$-double points 
lying on the curve $C\;:\; x=y=z=t$. 
Therefore, we obtain ${\rm rk} \, {\rm Cl}(X) = {\rm rk}\, 
{\rm Cl}(Y) = 9 - {\rm dim}\,T = 5$. 
Moreover, there exists exactly  one crepant divisor (over $C$) in a resolution 
of singularities of $X$. So the predicted power of $\log B$ in the asymptotic 
formula for $N(U, B)$ is $({\rm rk}\, {\rm Cl}(X) -1) + 1 = 5$. 
}
\end{exam}

\subsection{Cubic $xyz=u^3$}

We consider the singular cubic surface $X\subset \P^3$ over $\Q$
given by the homogeneous equation $xyz=u^3$. 
This is a toric variety, an 
equivariant compactification of the torus 
$$
T = {\rm Spec}
{\bf Q}[x,y,z]/(xyz - 1)
$$ 
given by the condition  $u\neq 0$. We can fix  an 
isomorphism $T \cong {\bf G}_m^2$ by choosing  $\{ x, y\} $ as 
a basis of the group of algebraic characters of $T$. 
Consider the problem of the computation of 
the asymptotic of 
$$
N(T,B) = \mbox{\rm Card}\{ (x,y) \in (\Q^*)^2\,\,:\,\,
 H(x,y) \le B\,\,\}
$$ 
for $B\ra \infty$, 
where 
\[ H(x,y) = \prod_{v \in {\rm Val}({\bf Q})} 
\max \{\|x\|_v,\|y\|_v,\|(xy)^{-1}\|_v,\|1\|_v \} \]
This problem is addressed in \cite{fouvry}. 
We would like to use this problem as a down-to-earth  illustration of  
 our general theory of 
height zeta functions of toric varieties.   
First of all we note that the relation $\|x\|_v \|y\|_v \|(xy)^{-1}\|_v = 
\| 1 \|_v =1$ implies that 
$$ 
H_v(x,y): = 
 \max \{\|x\|_v,\|y\|_v,\|(xy)^{-1}\|_v\}. 
$$  
Since $l_1 =\log  \|x\|_v$, $l_2 = \log \|y\|_v$, $l_3 = \log 
\|(xy)^{-1}\|_v$ are linear functions 
on the logarithmic space $N_{{\bf R},v} \cong {\bf R}^2$:  
$$
N_{{\bf R},v} =  
\left\{ \begin{array}{ll}  T({\bf Q}_v) /T({\cal O}_v) \otimes_{\bf Z}
{\bf R},   
& \mbox{\rm if $v = p \in {\rm Spec}\,{\bf Z}$} \\ 
 T({\bf Q}_v) /T({\cal O}_v),
 &   \mbox{\rm if $v = \infty$,}
\end{array} \right. 
$$ 
we can consider $\log h_v(x,y)$ as a piecewise linear function on 
 $N_{{\bf R},v}$. Let $e_1 = (-2,1)$, $e_2 = (1,-2)$ and 
$e_3 = (1,1)$ be  lattice vectors in ${\bf Z}^2 \subset {\bf R}^2$. 
We define the following $3$ convex cones in ${\bf R}^2$: 
\[ \sigma_1 = {\bf R}_{\geq 0} e_2 + {\bf R}_{\geq 0} e_3,  \] 
\[ \sigma_2 = {\bf R}_{\geq 0} e_1 + {\bf R}_{\geq 0} e_3,  \] 
\[ \sigma_3 = {\bf R}_{\geq 0} e_1 + {\bf R}_{\geq 0} e_2.  \]
Then $N_{{\bf R},v} = \bigcup_{i =1}^3 \sigma_i$ and the restriction 
of $\log H_v(x,y)$ to $\sigma_i$ coincides with the linear function $l_i$. 
Let 
\[  A := \bigoplus_{v \in {\rm Val}({\bf Q}} 
T({\bf Q}_v) /T({\cal O}_v) \]
be  the logarithmic adelic group. 
In order to compute the height zeta function
\[ Z(s) = \sum_{(x,y) \in {({\bf Q}^*)^2 = T({\bf Q})}} H(x,y)^{-s} \]
we use the natural homomorphism $Log$ of $T({\bf Q})$ 
to $A$. Denote by $B$ the subgroup $Log(T({\bf Q})) \subset A$.   
We remark that  the  kernel of $Log$ consists of $4$ elements of 
finite order in $({\bf Q}^*)^2$ and the quotient 
$A/B$ is isomorphic to ${\bf R}^2$. Moreover the functions 
$H_v(x,y)^{-s}$ on each  $T({\bf Q}_v) /T({\cal O}_v)$ define a natural  
extension of  $h(x,y)^{-s}$ to a function on $A$. So we obtain: 
\begin{equation}
 Z(s) = 4 \sum_{ b \in B} \prod_{v \in {\rm Val}({\bf Q})} H_v(b_v)^{-s} 
\label{poiss1}
\end{equation}

The main idea of our proof in 
\cite{BaTschi1} is to apply the Poisson summation
formula on the group $A$  and to express the
height zeta function $Z(s)$ as an integral 
$$
Z(s)=\frac{4}{(2\pi )^2}
\int_{{\R}^2} \left(
 \prod_{p} Q_p(s,i{\bf m})\cdot Q_{\infty}(s,i{\bf m}) \right) 
{\bf d}{\bf m},
$$
where ${\bf m} =(m_1, m_2) \in {\bf R}^2$, ${\bf dm} = dm_1dm_2$,  
$$
Q_{p}(s,i{\bf m}) = \sum_{(b_{1,p},b_{2,p}) \in {\bf Z}^2} 
h_p(b_p)^{-s} p^{ i<b, {\bf m}>}, 
$$
and
$$
Q_{\infty}(s,i{\bf m}) = \int_{{\bf R}^2} H_{\infty}(b)^{-s} 
\exp ( i<b, {\bf m}> ).   
$$

An  exact computation of $Q_{p}(s,i{\bf m})$ and 
$Q_{\infty}(s,i{\bf m})$ can be obtained by a subdivision 
of each of the  cones $\sigma_1, \sigma_2, \sigma_3$ into 
a union of $3$ subcones generated by a basis 
of the lattice ${\bf Z}^2 \subset {\bf R}^2$. 
(From the viewpoint of toric geometry this means that we 
reduce the 
counting problem for rational points in a torus 
with respect to a singular compactification 
to  
a counting problem for rational points in a torus 
with respect to the minimal resolution of 
singularities of this compactification). 
 This calculation is done
in \cite{BaTschi1}, Section 2 for arbitrary smooth toric varieties. 
What remains is the analytic continuation of the integral 
and the identification of the constant at the leading pole. 
For this it is necessary to work
on the whole complexified space $PL(\S)_{\C}$ and to invoke the
technical theorems in Section 6 in \cite{BaTschi2}. 
Applying the main theorem 
of \cite{BaTschi2}, we  obtain 
$$
N(T,B)= 
\frac{
\gamma_{{\cal K}^{-1}}(X)
 \delta_{{\cal K}^{-1}}(X) \tau_{{\cal K}^{-1}}(X)}{6!} 
B (\log B)^6(1+o(1))
$$
for $B\ra \infty$. The constants are as follows: 
$ \gamma_{{\cal K}^{-1}} (X)=1/36$, $\delta(X)=1$ 
and  $\tau_{{\cal K}^{-1}}(X)=\tau_{{\cal K}^{-1}}(X)_{\infty}
\prod_p\tau_{{\cal K}^{-1}}(X)_p$
where 
$$
\tau_{{\cal K}^{-1}}(X)_p= 
\left(1 + \frac{7}{p} + \frac{1}{p^2} \right)\cdot (1-1/p)^{7}
$$
for all primes $p$ and $\tau_{{\cal K}^{-1}}(X)_{\infty}=9\cdot 4$. 
Similar statements hold over any number field.
One can compute the constant
$\gamma_{{\cal K}^{-1}}(X)$ by observing that
$\Lambda_{\rm eff}^*(X)$ (the dual cone to the cone of
effective divisors) is a union of two
simplicial cones.

\bigskip

\noindent
\small  Mathematisches Institut, Universit\"at T\"ubingen   \\
\small  Auf der Morgenstelle 10,  72076  T\"ubingen, Germany  \\
\small  e-mail: batyrev@bastau.mathematik.uni-tuebingen.de \\
\small and \\
\small Dept. of Mathematics, U.I.C.\\
\small Chicago, (IL) 60607-7045,  U.S.A.  \\
\small e-mail: yuri@math.uic.edu

\end{document}